\documentclass[final]{siamart1116}%
\usepackage{amsmath}
\usepackage{amsfonts}
\usepackage{amssymb}
\usepackage{graphicx}%
\usepackage{lipsum}
\usepackage{amsfonts}
\usepackage{graphicx}
\usepackage{epstopdf}
\usepackage{algorithmic}
\usepackage{tikz}

\ifpdf
  \DeclareGraphicsExtensions{.eps,.pdf,.png,.jpg}
\else
  \DeclareGraphicsExtensions{.eps}
\fi

\usepackage{color}
\usepackage{listings}
\usepackage{hyperref}
\usepackage{listings}
\usepackage[framemethod=tikz]{mdframed}

\newcommand{\bu}{\mathbf{u}}
\newcommand{\bU}{\mathbf{U}}
\newcommand{\bt}{\mathbf{t}}

\newcommand{\bv}{\mathbf{v}}
\newcommand{\bV}{\mathbf{V}}

\newcommand{\ba}{\mathbf{a}}

\newcommand{\bb}{\mathbf{b}}

\newcommand{\bff}{\mathbf{f}}
\newcommand{\bg}{\mathbf{g}}

\newcommand{\bPhi}{\mathbf{\Phi}}

\newcommand{\bx}{\mathbf{x}}
\newcommand{\bX}{\mathbf{X}}

\newcommand{\bz}{\mathbf{z}}

\makeatletter
\newcommand*{\rom}[1]{\expandafter\romannumeral #1}
\makeatother

\numberwithin{theorem}{section}

\newcommand{\TheTitle}{$O(N)$ Hierarchical algorithm for computing the expectations of 
truncated multi-variate normal distributions in $N$ dimensions}
\newcommand{\TheAuthors}{J. Huang, F. Fang, G. Turkiyyah, J. Cao, M.G. Genton, and D.E. Keyes }

\headers{Hierarchical algorithm for TMVN expectations}{\TheAuthors}

\title{{\TheTitle}\thanks{Submitted to the editors .}}

\author{
	Jingfang Huang, Fuhui Fang \\ \vspace{-0.05in} {\it \scriptsize University of North Carolina at Chapel Hill} \\  
  George Turkiyyah \\ \vspace{-0.05in} {\it \scriptsize American University of Beirut} \\
  Jian Cao, Marc G. Genton, David E. Keyes \\ \vspace{-0.1in} {\it \scriptsize King Abdullah University of Science and Technology} 
}

\usepackage{amsopn}

\ifpdf
\hypersetup{
  pdftitle={\TheTitle},
  pdfauthor={\TheAuthors}
}
\fi

\begin{document}

\maketitle

\begin{abstract}
In this paper, we study the $N$-dimensional integral 
$\phi(\ba, \bb; A) = \int_{\ba}^{\bb} H(\bx) f(\bx | A) \text{d} \bx$ representing 
the expectation of a function $H(\bX)$ where $f(\bx | A)$ is the truncated multi-variate 
normal (TMVN) distribution with zero mean, $\bx$ is the vector of integration 
variables for the $N$-dimensional random vector $\bX$, $A$ is the inverse of the 
covariance matrix $\Sigma$, and $\ba$ and 
$\bb$ are constant vectors. We present a new hierarchical algorithm which can 
evaluate $\phi(\ba, \bb; A)$ using asymptotically optimal $O(N)$ operations 
when $A$ has ``low-rank" blocks with ``low-dimensional" features and 
$H(\bx)$ is ``low-rank". 
We demonstrate the divide-and-conquer idea when $A$ is a symmetric positive 
definite tridiagonal matrix, and present the necessary building blocks 
and rigorous potential theory based algorithm analysis when $A$ is given 
by the {\it exponential covariance model}. Numerical results are presented 
to demonstrate the algorithm accuracy and efficiency for these two cases. 
We also briefly discuss how the algorithm can be generalized to a wider 
class of covariance models and its limitations. 
\end{abstract}



\section{Introduction}
In this paper, we study the efficient computation of the expectation of function
$H(\bX)$ given by
\begin{eqnarray}
\label{eq:gaussian}
\phi(\ba, \bb; A) &= & \int_{\ba}^{\bb} H(\bx) f(\bx | A) \text{d}\bx \nonumber \\
& = & 
  \int_{a_1}^{b_1} \cdots \int_{a_N}^{b_N} H(\bx) |\Sigma|^{-1/2} (2 \pi)^{-N/2} 
  e^{ -\frac12 \bx^T A \bx } \text{d}x_N \cdots \text{d}x_1, 
\end{eqnarray}
where the $N$-dimensional random vector $\bX=(X_1, \ldots, X_N)^T$ follows 
the truncated multivariate normal distribution (TMVN), $f(\bx|A)$ is the 
$N$-dimensional multivariate Gaussian probability density function with zero mean 
and covariance matrix $\Sigma$, $A$ is the inverse of the symmetric positive 
definite (SPD) $N\times N$ covariance matrix $\Sigma$,  $\bx$ is the integration 
variable, and the integration limits are $\ba=(a_1, \ldots, a_N)^T$ and 
$\bb=(b_1,\ldots, b_N)^T$ which form a hyper-rectangle 
in $\Bbb{R}^N$. The efficient computation of $\phi$ is very important for many 
applications, including those in spatial and temporal statistics and in the study
of other high dimensional random data sets where the Gaussian distribution is commonly 
used, see \cite{arellano2006unification,arellano2006unified,arellano2008exact,azzalini2013skew,
castruccio2016high,dombry2017full,genton2004skew,genton2018hierarchical,stephenson2005exploiting} 
and references therein. Note that $|\Sigma|^{-1/2}(2\pi)^{-N/2}$ is a constant. When
$\Sigma$ has low-rank properties, $|\Sigma|^{-1/2}$ can be evaluated efficiently using 
existing low-rank linear algebra techniques \cite{greengard2009fast,hackbusch1999sparse,
hackbusch2000sparse,ho2012fast}. we ignore this term to simplify our discussions in this paper. 

Due to the ``curse of dimensionality", direct evaluation of the $N$-dimensional integral
using standard quadrature rules is computationally demanding (and impossible for many 
settings using today's supercomputers), and most existing schemes either scale poorly when
the dimension $N$ increases or rely on the Monte Carlo methods for very high dimensional cases
\cite{botev2017normal,craig2008new,genz1992numerical,genz2009computation,genz2014multivariate,
geweke1991efficient,keane199320,hajivassiliou1996simulation,meyer2010recursive,miwa2003evaluation,
phinikettos2011fast,ridgway2016computation}. A good review of
existing techniques can be found in \cite{genz2009computation}. The purpose of this
paper is to show that when there exist special structures in $H$ and $A$ (or equivalently
in $\Sigma$), fast direct evaluation of the $N$-dimensional integral becomes possible. 
In particular, when the function $H(\bx)$ is ``low-rank" and the matrix $A$ has 
hierarchical low-rank blocks with ``low-dimensional" singular vectors in their 
singular value decompositions, 
asymptotically optimal $O(N)$ {\it hierarchical} algorithms can be developed, by 
compressing these compact features and efficiently processing them ``locally" on 
a hierarchical tree structure. We leave the mathematical rigorous definitions of the
``low-rank" and ``low-dimensional" concepts to later discussions, but only
mention that such ``compact" structures exist in many important applications. For
example, when the high dimensional data can be properly clustered, e.g., by using
their spatial or temporal locations and relative distance or pseudo-distance, 
the matrix blocks describing the ``interactions" between different clusters are 
often low-rank as revealed by the principal component analysis (PCA). 

This paper presents the algorithm analysis and implementation details for two representative
matrices: (a) when $A$ is a tridiagonal SPD matrix; and (b) when $A$ has the same form as
the covariance matrix in the exponential covariance model in one dimensional setting. 
In case (b) when $A$ is the exponential covariance matrix, the original covariance 
matrix $\Sigma=A^{-1}$ is approximately a tridiagonal system. In the numerical algorithm 
for both cases, a downward pass is first performed on a hierarchical
tree structure, by introducing a $t$-variable to divide the parent problem (involving a 
function with no more than $2$ ``effective" variables) into two child problems, each involving 
a function with no more than $2$ ``effective" variables. The relation coefficients between 
the parent's
effective variables, new $t$-variable, and children's effective variables are computed and 
stored for each tree node. At the leaf level, the one dimensional integral which only involves 
one $x_j$ variable is evaluated either analytically or numerically, and then
approximated numerically by a global Fourier series representation. An upward pass is then 
performed, to recursively forming the approximating Fourier series of the parent's $2$-effective 
variable function using those from its two children. The function value $\phi$
in Eq.~(\ref{eq:gaussian}) is simply given by the constant function (with two ``null variables")
at the root level of the tree structure.
The presented hierarchical algorithms share many similar features as many existing
fast hierarchical algorithms in scientific computing, including the classical fast Fourier 
transform (FFT) \cite{cooley1965algorithm}, multigrid method (MG) 
\cite{brandt1977multi,hackbusch2013multi}, fast multipole method (FMM) 
\cite{greengard1987fast,greengard1997new}, and the fast direct solvers (FDS) and 
hierarchical matrix ($\mathcal{H}$-matrix) algorithms 
\cite{greengard2009fast,ho2012fast,hackbusch1999sparse,hackbusch2000sparse}. 

This paper is organized as follows. In Sec.~\ref{sec:concept}, we introduce the mathematical
definitions of the ``low-rank" and ``low-dimensional" concepts. In Sec.~\ref{sec:tridiag}, 
we present the details of a hierarchical algorithm for computing $\phi$ when $A$ is a 
tridiagonal matrix. In Sec.~\ref{sec:rankone}, we show how the algorithm can be generalized 
to the case when $A$ has the same form as the exponential covariance matrix in one 
dimensional setting, and present the rigorous analysis using potential
theory from ordinary and partial differential equation analysis, as the exponential covariance 
model in one dimension is closely related with the Green's function and integral equation
solutions of the boundary value ordinary differential equation $u(x)-u''(x)=f(x)$.
Numerical results are presented to demonstrate the accuracy, stability, and $O(N)$ complexity 
of the new hierarchical algorithm for both cases. 
In Sec.~\ref{sec:generalize}, we discuss how the algorithm can be generalized to more complicated 
cases as well as its limitations. In particular, our current algorithm implementation relies 
heavily on existing numerical tools and software packages for accurately processing multi-variable functions 
(e.g., high dimensional non-uniform FFT or sparse grid techniques). Many of these tools are
unfortunately still unavailable even when the number of independent variables is approximately $5 \sim 20$.
Finally in Sec.~\ref{sec:conclusion}, we summarize our results.

\section{Low-rank and Low-dimensional Properties}
\label{sec:concept}
Our algorithm can be applied to a function $H(\bx)$ with the following structure,
\begin{equation}
\label{eqn:H-lowrank}
H(\bx) = \sum_{p=1}^P u_{p,1}(x_1) u_{p,2}(x_2) \cdots u_{p,N}(x_N)=\sum_{p=1}^P \prod_{k=1}^N u_{p, k}(x_k),
\end{equation}  
where $P$ is assumed to be a small constant independent of $N$, and each function $u_{p, k}$ 
is a single variable function, not necessarily a continuous function. As the separation
of variables 
$$H(x,y)=\sum_{p=1}^P u_p(x) v_p(y) $$ 
can be considered as the non-orthogonalized function version of the singular value 
matrix decomposition
$$H_{m \times n} = U_{m \times P} \Lambda_{P \times P} V^T_{P \times n},$$
we refer to a function $H$ with a representation in Eq.~(\ref{eqn:H-lowrank}) as a 
{\bf low-rank (rank-$P$) function}.
Plugging Eq.~(\ref{eqn:H-lowrank}) into Eq.~(\ref{eq:gaussian}), the original problem
of evaluating $\phi$ now becomes the evaluations of $P$ integrals, each has the form
\begin{equation}
\label{eq:GaussNew}
\phi_p(\ba, \bb; A) = 
  \int_{a_1}^{b_1} \cdots \int_{a_N}^{b_N} \prod_{k=1}^N u_{p, k}(x_k) 
  \exp \left(-\frac12 \bx^T A \bx \right) \text{d}x_N \cdots \text{d}x_1.
\end{equation}
We focus on $\phi_p$ in the following discussions, and simply denote $\phi_p$ as 
$\phi$.

\begin{figure}[htbp]
  \centering
  \begin{tikzpicture} 
\fill[color=blue!40] (0,2) rectangle (2,4)
(2,0) rectangle (4,2);
\fill[color=green!40] (0,0) rectangle (2,2)
(2,2) rectangle (4,4); 
\draw[step=2] (0,0) grid (4,4);
\node at (1,3) {$A_{1,1}$};
\node at (3,3) {$A_{1,2}$};
\node at (1,1) {$A_{1,3}$};
\node at (3,1) {$A_{1,4}$};
\fill[color=blue!40] (6,3) rectangle (7,4)
(7,2) rectangle (8,3)
(8,2) rectangle (9,1)
(9,1) rectangle (10,0);
\fill[color=green!40] (6,0) to (6,3) to (7,3) to (7,2)
to (8,2) to (8,1) to (9,1) to (9,0) to (6,0);
\fill[color=green!40] (7,4) to (10,4) to (10,1)
to (9, 1) to (9, 2) to (8, 2) to (8, 3) to (7,3) to (7,4);  
\draw[step=2] (6,0) grid (10,4);
\draw[step=1] (6,2) grid (8, 4);
\draw[step=1] (8, 0) grid (10, 2);
\node at (6.5,3.5) {$A_{2,1}$};
\node at (7.5,3.5) {$A_{2,2}$};
\node at (6.5,2.5) {$A_{2,3}$};
\node at (7.5,2.5) {$A_{2,4}$};
\node at (8.5,1.5) {$A_{2,5}$};
\node at (9.5,1.5) {$A_{2,6}$};
\node at (8.5,0.5) {$A_{2,7}$};
\node at (9.5,0.5) {$A_{2,8}$};
\end{tikzpicture} 
  \caption{$\mathcal{H}$-matrix after $2$ (left) and $3$ (right) divisions, 
	 with low-rank off-diagonal blocks (green).}
  \label{fig:hmatrix}
\end{figure}
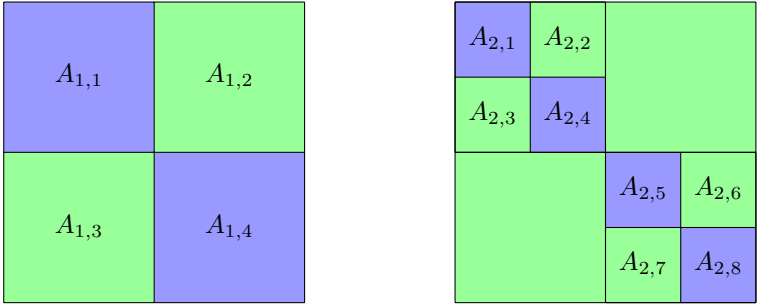

For the inverse $A$ of the covariance matrix $\Sigma$, we assume it belongs 
to a class of hierarchical matrices ($\mathcal{H}$-matrices) 
\cite{hackbusch1999sparse,hackbusch2000sparse} with low-rank off-diagonal
blocks. A sample Hierarchical matrix after $2$ (left) and $3$ (right) divisions is 
demonstrated in Fig.~\ref{fig:hmatrix}, where the blue square block represents 
the self-correlation within each cluster of random variables $X_i$, and the 
green block shows the correlation
between two different clusters. We define a cluster in the original domain as 
a set of indices of the column vectors, and a cluster in the target space as a
set of indices of the row vectors. The correlation between the cluster $J$ of 
the original domain and cluster $K$ of the target space is described by the 
matrix block formed by only extracting the $K$-entries from the $J$-columns.  
We consider $\mathcal{H}$-matrices with {\bf low-rank off-diagonal blocks}, by assuming 
that the ranks of all the off-diagonal blocks are bounded by a constant $P$, 
which is independent of the block matrix size. We use $P$ to represents the dimension
of a subspace or the rank of a matrix in this paper, and the rank $P$ of the off-diagonal 
blocks can be different from the rank $P$ in Eq.~(\ref{eqn:H-lowrank}). 
In numerical linear algebra, ``low-rank off-diagonal block" means that the 
off-diagonal block $A_{i,j}$ of size $n \times n$ has the following singular value 
decomposition
$$A_{i,j} = U_{i,j} \Lambda_{i,j} V^T_{i,j},$$
where $U$ and $V$ are of size $n \times P$ and respectively contain the orthonormal 
vectors in the target space and original domain, and $\Lambda$ is a size 
$P \times P$ diagonal matrix with ordered and non-negative diagonal entries. 
As the random variables $\{X_i, i=1,\ldots, N \}$ are clustered hierarchically, we index the block matrices 
$A_{i,j}$ differently from those commonly used in matrix theory to emphasize this 
hierarchical structure in the $\mathcal{H}$-matrix, where $i$ represents the level of the matrix 
block, and $j$ is its index in that particular level. The original matrix $A$ is defined
as the level $0$ matrix. After the $1^{st}$ division, the $4$ matrix blocks are 
indexed $(1,1)$, $(1,2)$, $(1,3)$, and $(1,4)$. The diagonal matrix blocks will be further
divided and the off-diagonal matrices become leaf nodes to form an {\bf adaptive quad-tree
structure}. In the left of Fig.~\ref{fig:hmatrix}, the matrix $A_{1,2}$ denotes the 
second matrix block at level $1$, representing the correlations between the second 
cluster in the original domain and first cluster in the target space. 
For the covariance matrix, as the target space and original domain are the one 
and the same, the indices of the random variables $\bX$ (and integration variables $\bx$)
will be used to cluster the indices of both the target space and original domain, to form a 
{\bf uniform binary tree structure}. In the following, we focus on the integration variables 
$\{x_i, i=1,\ldots, N\}$, which are referred to as the $x$-variables.

Next, we consider the ``low-dimensional" concept, by studying a function with $M$ 
$t$-variables $t_1$, $t_2$, $\ldots$, $t_M$ of the form
$$F(t_1 \bu_1+t_2 \bu_2+ \cdots +t_M \bu_M).$$
When the dimension $P$ of the vector space $span\{\bu_1, \bu_2, \ldots, \bu_M\}$ is 
much less than $M$, $P << M$, we say $F$ is a ``low-dimensional" function.
Assuming the basis for the vector space $span\{\bu_1, \bu_2, \ldots, \bu_M\}$
is given by $\{\bv_1, \bv_2, \ldots, \bv_P\}$, the function $F$ can be considered as
an ``effective" $M$ variable function, where the new $w$-variables 
$\{ w_1, w_2,\ldots, w_P\}$ are combinations of the $t$-variables and
satisfy the relation 
$$w_1 \bv_1+w_2 \bv_2+ \cdots +w_P \bv_P = t_1 \bu_1+t_2 \bu_2+ \cdots +t_M \bu_M.$$

The low-rank and low-dimensional structures exist in many practical systems. The
well studied low-rank concept measures the rank of a matrix block and is closely 
related with the principal component analysis in statistics and singular value 
decomposition (SVD) in numerical linear algebra. When the data can be clustered,
the covariance matrix block describing the relations between two different clusters
is often low-rank, and both the storage of such a matrix block and related operations 
can be reduced significantly using today's low-rank linear algebra techniques.
The {\it low-dimensional} property in this paper considers the special structures in the
singular vectors of the SVD decomposition of the low-rank off-diagonal blocks.
Consider two clusters of the $x$-variables and the space formed by extracting 
all the corresponding sub-vectors describing the relations of these two clusters 
from the singular vectors in the SVD decompositions of all the off-diagonal matrices.
When the covariance matrix is defined by a covariance function using the spatial or 
temporal locations (or pseudo-locations) $z_i$ and $z_j$ of the corresponding 
random variables $X_i$ and $X_j$, the covariance function is often ``smooth" 
and only contains ``low-frequency" information when $i \neq j$, it 
can be well approximated by a few terms of truncated Taylor expansion 
(or other basis functions) when a separation of variables is performed on the 
covariance function determined by the two location variables $z_i$ and $z_j$.  
In this case, all the singular vectors are the discretized versions of the 
polynomial basis functions at locations corresponding to the cluster index 
sets. The dimension of the space formed by these singular vectors is 
therefore determined by the highest degree of the polynomial basis functions. 
When the $\bu_i$ vectors are extracted from these singular vectors,
the function $F(t_1 \bu_1+t_2 \bu_2+ \cdots +t_M \bu_M)$ will be low-dimensional
and the number of effective variables is also determined by the highest degree
of the polynomial basis functions. The special structures in the singular 
vectors were also used in \cite{ho2012fast,ho2016hierarchical,l2016hierarchical}.
The low-rank and low-dimensional concepts will be further studied in the next two sections.
 
\section{Case I: Tridiagonal System}
\label{sec:tridiag}
We demonstrate the basic ideas of the hierarchical algorithm by studying a simple
tridiagonal system 
\begin{equation}
A = \left[ \begin{array}{rrrrrrr}
4 & -2 & 0 & . & . & . & 0\\
-2 & 4 & -2 & 0 & . & . & 0\\
0 & -2 & 4 & -2 & 0 & . & 0\\
.& . & . & . & . & . & . \\
. & . & . & . & . & . & . \\
. & . & . & . & . & . & . \\
0 & . & . & . & 0 & -2 & 4\end{array} \right]_{N \times N}.
\end{equation}
We assume $N=2^L$ and first consider a constant function $H(\bx)$ to 
simplify the notations and discussions. The algorithm for more general 
low-rank $H(\bx)$ in Eq.~(\ref{eq:GaussNew}) only requires a slight 
change in the code for the leaf nodes, which will become clear after 
we present the algorithm details for the simplified integration problem
\begin{eqnarray}
\label{eq:tridiag}
\phi(\ba, \bb; A) & = &  
  \int_{a_1}^{b_1} \cdots \int_{a_N}^{b_N} 
  e^{ - \frac12 \bx^T A \bx } \text{d}x_N \cdots \text{d}x_1 \nonumber \\
  & = &
  \int_{\ba}^{\bb} 
  e^{ \left(-2x_1^2 + 2 x_1x_2 -2 x_2^2 +  \cdots - 2 x_k^2 + {\color{red} 2 x_k x_{k+1}} -2 x_{k+1}^2 +
  \cdots - 2 x_{N-1}^2 + 2 x_{N-1} x_{N} -2 x_N^2 \right)} \text{d}\bx,
\end{eqnarray}
where $k=2^{L-1}=N/2$.  
The tridiagonal matrix is a very special $\mathcal{H}$-matrix, where each off-diagonal
matrix block only contains one non-zero number either at the lower-left or 
upper-right corner of the matrix block and is rank $1$. The singular
vectors are either $\bu_i=[1,0,0,\ldots,0]^T$ or $\bu_i=[0,0,\ldots, 0,1]^T$. For
any given cluster of indices, the number of effective variables in 
$t_1 \bu_1 + t_2 \bu_2 + \cdots + t_M \bu_M$ is therefore no more than 2, and
the only non-zero numbers are located either at the first or the last entry 
in the singular vectors $\bu_1$, $\bu_2$, $\ldots$, $\bu_M$. 

\subsection{Divide and Conquer on a Hierarchical Tree}
Note that the $x$-variables $[x_1, \ldots, x_k]$ and $[x_{k+1}, \ldots, x_N ]$ are coupled in the
integrand only through {\bf one} term $2 x_k x_{k+1}$. If this ``weak coupling" term had not 
been there, then we would have two completely decoupled ``child problems", 
and the integral could be evaluated as
\begin{displaymath}
\begin{array}{l}
 \int_{\ba}^{\bb} 
  e^{ \left(-2x_1^2 + 2 x_1x_2 -2 x_2^2 +  \cdots - 2 x_k^2 -2 x_{k+1}^2 + \cdots
 - 2 x_{N-1}^2 + 2 x_{N-1} x_{N} -2 x_N^2 \right)} \text{d}\bx \\ 
=  \left( \int_{a_1}^{b_1} \cdots \int_{a_k}^{b_k}  
  e^{ \left(-2x_1^2 + 2 x_1x_2 -2 x_2^2 +  \cdots - 2 x_k^2 \right)} \text{d}x_k \cdots \text{d}x_1 \right) 
  \cdot \\ \left( \int_{a_{k+1}}^{b_{k+1}} \cdots \int_{a_N}^{b_N}  
  e^{ \left( -2 x_{k+1}^2 + \cdots
 - 2 x_{N-1}^2 + 2 x_{N-1} x_{N} -2 x_N^2 \right) } \text{d}x_N \cdots \text{d}x_{k+1} \right).
\end{array}
\end{displaymath}
If the same assumptions could be made to each ``child problem", then the high-dimensional integral
would become the product of $N$ one-dimensional integrals. 

A convenient tool to decouple the $x$-variables in order to have two child problems is to use the
Fourier transform formula for the Gaussian distribution as 
\begin{equation}
\label{eq:Fourier}
  e^{-(x-y)^2} = c \int_{-\infty}^{\infty} e^{2 i t (x-y)} e^{- t^2} \text{d}t =  c \int_{-\infty}^{\infty} 
{\color{red} e^{2 i t x} e^{ - 2 i t y} } e^{- t^2} \text{d}t
\end{equation}
where $i=\sqrt{-1}$ and $c=\frac{1}{\sqrt{\pi}}$. Note that the variables $x$ and $y$ are decoupled 
on the right hand side. Completing the square in Eq.~(\ref{eq:tridiag}) and applying 
the formula in Eq.~(\ref{eq:Fourier}) to the resulting $(x_k-x_{k+1})^2$ term (in red) in the
integral, we get
\begin{displaymath}
\begin{array}{l}
 \phi(\ba, \bb; A)=\int_{\ba}^{\bb} 
  e^{ \left(-2x_1^2 + 2 x_1x_2 -2 x_2^2+  \cdots {\color{blue} - 2 x_k^2 + 2 x_k x_{k+1} -2 x_{k+1}^2 } 
  	 +
  	 \cdots
 - 2 x_{N-1}^2 + 2 x_{N-1} x_{N} -2 x_N^2 \right)} \text{d}\bx \\ 
 = \int_{\ba}^{\bb} e^{ \left(-2x_1^2 + 2 x_1x_2 -2 x_2^2 +  \cdots 
       {\color{blue} -  x_k^2 {\color{red} - (x_k - x_{k+1})^2} - x_{k+1}^2 } + \cdots
 - 2 x_{N-1}^2 + 2 x_{N-1} x_{N} -2 x_N^2 \right)} \text{d}\bx \\ 
 = c \int_{-\infty}^{\infty} e^{- t^2} h_{1,1}(t) h_{1,2}(t) \text{d}t 
\end{array}
\end{displaymath}
where $h_{1,1}(t)$ and $h_{1,2}(t)$ are both single $t$-variable functions given by
\begin{displaymath}
\begin{array}{l}
 h_{1,1}(t) = \int_{a_1}^{b_1} \cdots \int_{a_k}^{b_k}  
  e^{\left(-2x_1^2 + 2 x_1x_2 -2 x_2^2 + \cdots  {\color{blue} - x_k^2 + 2 i x_k t} \right)} \text{d}x_k \cdots \text{d}x_1,\\
 h_{1,2}(t) = \int_{a_{k+1}}^{b_{k+1}} \cdots \int_{a_N}^{b_N}  
  e^{\left( {\color{blue} - 2 i x_{k+1} t - x_{k+1}^2} + \cdots
 - 2 x_{N-1}^2 + 2 x_{N-1} x_{N} -2 x_N^2 \right) } \text{d}x_N \cdots \text{d}x_{k+1}.
\end{array}
\end{displaymath}
Note that the $x$-variables in the original problem (associated with a root node at level $0$ of a 
binary tree structure) are divided into two subsets of the same size, each set is associated with 
a ``child node" and a single $t$-variable function $h_{1,k}(t)$, $k=1$ or $k=2$. 

By introducing two new $t$-variables $t_{1,1}$ and $t_{1,2}$ for the functions 
$h_{1,1}$ and $h_{1,2}$, respectively, the same technique can be applied to decouple the 
$x$-variables $[x_1, \ldots, x_{\frac{k}{2}}]$ and $[x_{\frac{k}{2}+1}, \ldots, x_k ]$ in 
$h_{1,1}(t)$ and the $x$-variables  $[x_{k+1}, \ldots, x_{\frac{3k}{2}}]$ 
and $[x_{\frac{3k}{2}+1}, \ldots, x_N ]$ in $h_{1,2}(t)$, to derive
\begin{eqnarray*} 
h_{1,1}(t)=c \int_{-\infty}^{\infty} e^{- t_{1,1}^2} h_{2,1}(t_{1,1}) h_{2,2}(t_{1,1}, t) \text{d}t_{1,1}  \\
h_{1,2}(t)=c \int_{-\infty}^{\infty} e^{- t_{1,2}^2} h_{2,3}(t,t_{1,2}) h_{2,4}(t_{1,2}) \text{d}t_{1,2}, 
\end{eqnarray*}
where 
\begin{displaymath}
\begin{array}{l}
 h_{2,1}(t_{1,1}) = c \int_{a_1}^{b_1} \cdots \int_{a_{\frac{k}{2}}}^{b_\frac{k}{2}}  
  e^{\left(-2x_1^2 + 2 x_1x_2 -2 x_2^2 + \cdots {- x_{\frac{k}{2}}^2 
	+ 2 i x_{\frac{k}{2} } t_{1,1}} \right)} \text{d}x_{\frac{k}{2}} \cdots \text{d}x_1,\\
 h_{2,2}(t_{1,1},t) = c \int_{a_{\frac{k}{2}+1}}^{b_{\frac{k}{2}+1}} \cdots \int_{a_k}^{b_k}  
  e^{\left(-2 i x_{\frac{k}{2}+1} t_{1,1} - x_{\frac{k}{2}+1}^2 +  \cdots - x_k^2 + 2 i x_k t \right)} 
   	\text{d}x_k \cdots \text{d}x_{\frac{k}{2}+1},\\
 h_{2,3}(t,t_{1,2}) = c \int_{a_{k+1}}^{b_{k+1}} \cdots \int_{a_{\frac{3k}{2}}}^{b_{\frac{3k}{2}}}  
  e^{\left( - 2 i x_{k+1} t - x_{k+1}^2 + \cdots
 - x_{\frac{3k}{2}}^2 + 2 i x_{\frac{3k}{2}} t_{1,2} \right) } \text{d}x_{\frac{3k}{2}} \cdots \text{d}x_{k+1}, \\
 h_{2,4}(t_{1,2}) = c \int_{a_{\frac{3k}{2}+1}}^{b_{\frac{3k}{2}+1}} \cdots \int_{a_N}^{b_N}  
  e^{\left( { - 2 i x_{\frac{3k}{2}+1} t_{1,2} - x_{\frac{3k}{2}+1}^2} +\cdots
  + 2 x_{N-1} x_{N} -2 x_N^2 \right) } \text{d}x_N \cdots \text{d}x_{\frac{3k}{2}+1},
\end{array}
\end{displaymath}

\begin{figure}[htbp]
  \centering
  \begin{tikzpicture}[
    level/.style={sibling distance = 2cm/#1, level distance = 2.5cm},
    grow=right, 
  ] 
  \node {$\phi(\ba,\bb;A)$}
  child { node [scale=0.9] {$h_{1,1}(t_0)$} 
    child { node [scale=0.8] {$h_{2,1}(t_{1,1})$} 
      child { node [scale=0.7] (31) {$h_{3,1}(t_{2,1})$}
      }
      child { node [scale=0.7] (32) {$h_{3,2}(t_{2,1}, t_{1,1})$}
      }
    }
    child { node [scale=0.8] {$h_{2,2}(t_{1,1}, t_0)$}
      child { node [scale=0.7] (33) {$h_{3,3}(t_{1,1}, t_{2,2})$}
      }
			child { node [scale=0.7] (34) {$h_{3,4}(t_{2,2}, t_0)$}
      }
    }                            
  }
  child { node  [scale=0.9] {$h_{1,2}(t_0)$}
    child { node [scale=0.8] {$h_{2,3}(t_0, t_{1,2})$}
      child { node [scale=0.7] (35) {$h_{3,5}(t_0, t_{2,3})$}
      }
      child { node [scale=0.7] (36) {$h_{3,6}(t_{2,3}, t_{1,2})$}
      }
    } 
    child { node [scale=0.8] {$h_{2,4}(t_{1,2})$} 
      child { node [scale=0.7] (37) {$h_{3,7}(t_{1,2}, t_{2,4})$} 
      }
      child { node [scale=0.7] (38) {$h_{3,8}(t_{2,4})$}
      }
    } 
	}; 
\node[right of=31, xshift=15pt, scale=0.7] {$=\int e^{-x_1^2+2it_{2,1}x_1} \text{d}x_1$};
\node[right of=32, xshift=25pt, scale=0.7] {$=\int e^{2ix_2(t_{1,1}-t_{2,1})}\text{d}x_2$};
\node[right of=33, xshift=25pt, scale=0.7] {$=\int e^{2ix_3(t_{2,2} - t_{1,1})} \text{d}x_3$};
\node[right of=34, xshift=20pt, scale=0.7] {$=\int e^{2ix_4(t_0-t_{2,2})}\text{d}x_4$};
\node[right of=35, xshift=20pt, scale=0.7] {$=\int e^{2ix_5(t_{2,3} - t_0)} \text{d}x_5$};
\node[right of=36, xshift=25pt, scale=0.7] {$=\int e^{2ix_6(t_{1,2}-t_{2,3})}\text{d}x_6$};
\node[right of=37, xshift=25pt, scale=0.7] {$=\int e^{2ix_7(t_{2,4} - t_{1,2})} \text{d}x_7$};
\node[right of=38, xshift=15pt, scale=0.7] {$=\int e^{-x_8^2 - 2ix_8t_{2,4}}\text{d}x_8$};
\end{tikzpicture}
  \caption{A three-level partition that decomposes the original $N$-dimensional ($N=8$)
    integral. }
  \label{fig:p8}
\end{figure}
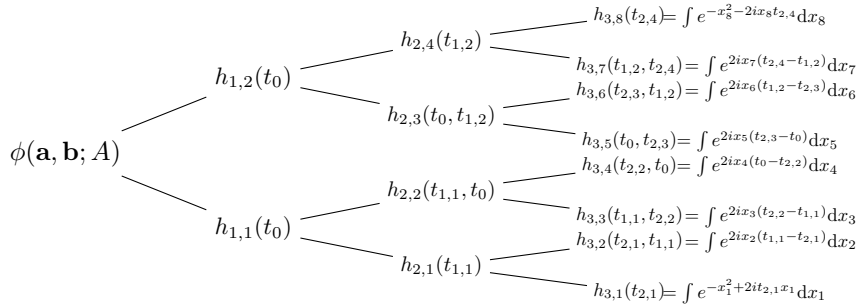

Repeating this procedure recursively on the hierarchical tree structure derived by recursively 
dividing the parent's $x$-variable set into two child subsets of the same size, a hierarchical
$h$-function $h_{l,k}$ will be defined for each tree node, where $\{l, k\}$ is the index of 
the tree node defined in the same way as that of the $x$-variable sets. 
One can show that for a parent node with index $p$,
its $h$-function $h_p(t_l,t_r)$ (with at most two $t$-variables $t_l$ and $t_r$) can be computed 
from the two child functions $h_{c_1}(t_l,t_m)$ and $h_{c_2}(t_m,t_r)$ (each with at most two 
$t$-variables) by integrating the $t$-variable $t_m$ used to decouple the parent problem using 
Eq.~(\ref{eq:Fourier}) as
\begin{equation}
\label{eq:L2L}
{
h_p(t_l,t_r)=c \int_{-\infty}^{\infty} e^{- t_m^2} h_{c_1}(t_l,t_m) h_{c_2}(t_m,t_r) \text{d}t_m.  }
\end{equation}
At the finest level when the $x$-variable set only contains one $x$-variable $x_j$, the two $t$-variable function
is given by 
$$ h_{leaf_{x_j}} (t_l,t_r) = c \int_{a_j}^{b_j} e^{\alpha x_j^2 -2 i x_j (t_l-t_r)} \text{d} x_j $$
where $\alpha=0$ for the interior nodes and $\alpha=-1$ for the two boundary nodes at the leaf level. For each
boundary node in the tree structure, its associated $h$-function only involves one $t$-variable
as the other becomes a null variable. 
In Fig.~\ref{fig:p8}, we show the detailed decoupling procedure and the functions $h_{j,k}$ when $N=8$,
where the first index $j$ of $t_{j,k}$ indicates the level at which the new $t$-variable is introduced, 
and the second index $k$ is its index at this level, ordered from bottom (left boundary of $x$-variables) 
to top (right boundary) in the figure. 

\vspace{0.1in}
{\noindent \bf Remark:} Each parent's $h$-function has no more than {\it two} 
$t$-variables, and it can be computed using the two children's $h$-functions, each with no 
more than {\it two} $t$-variables, as shown in Eq.~(\ref{eq:L2L}). Note that
the decoupling process is performed on a hierarchical binary tree structure, by introducing one
new $t$-variable and dividing parent's $x$-variable set into two children's subsets of the same 
size. As the depth of the tree is $O(\log N)$ so a total of $O(\log N)$ $t$-variables
will be introduced for each tree branch from the root to leaf level. However,
as the singular vectors are either $\bu_i=[1,0,0,\ldots,0]^T$ or $\bu_i=[0,0,\ldots, 0,1]^T$. For a 
tree node containing a particular set of $x$-variable indices from $x_{j+1}$ to $x_{j+k}$, 
there are at most two non-zero vectors in the vector set $\{ \bu_1, \bu_2, \ldots, \bu_M \}$, 
with the non-zero entry located either at the first or the last entry in one of the two 
non-zero singular vectors of size $k$. The number of effective variables in 
$t_1 \bu_1 + t_2 \bu_2 + \cdots + t_M \bu_M$ is therefore no more than 2, 
and 
$$[x_{j+1}, x_{j+2}, \ldots, x_{j+k}] \cdot \left( t_1 \bu_1 + t_2 \bu_2 + \cdots + t_M \bu_M \right)
  = -2 i x_{j+1} t_l + 2 i x_{j+k} t_r. $$
Therefore all the $h$-functions in the hierarchical tree structure have no more than
two effective variables and are ``low-dimensional" functions.

\subsection{Algorithm Details}
\begin{figure}[ht!]
\centering
\includegraphics[height=0.3\textwidth,width = 0.8\textwidth]{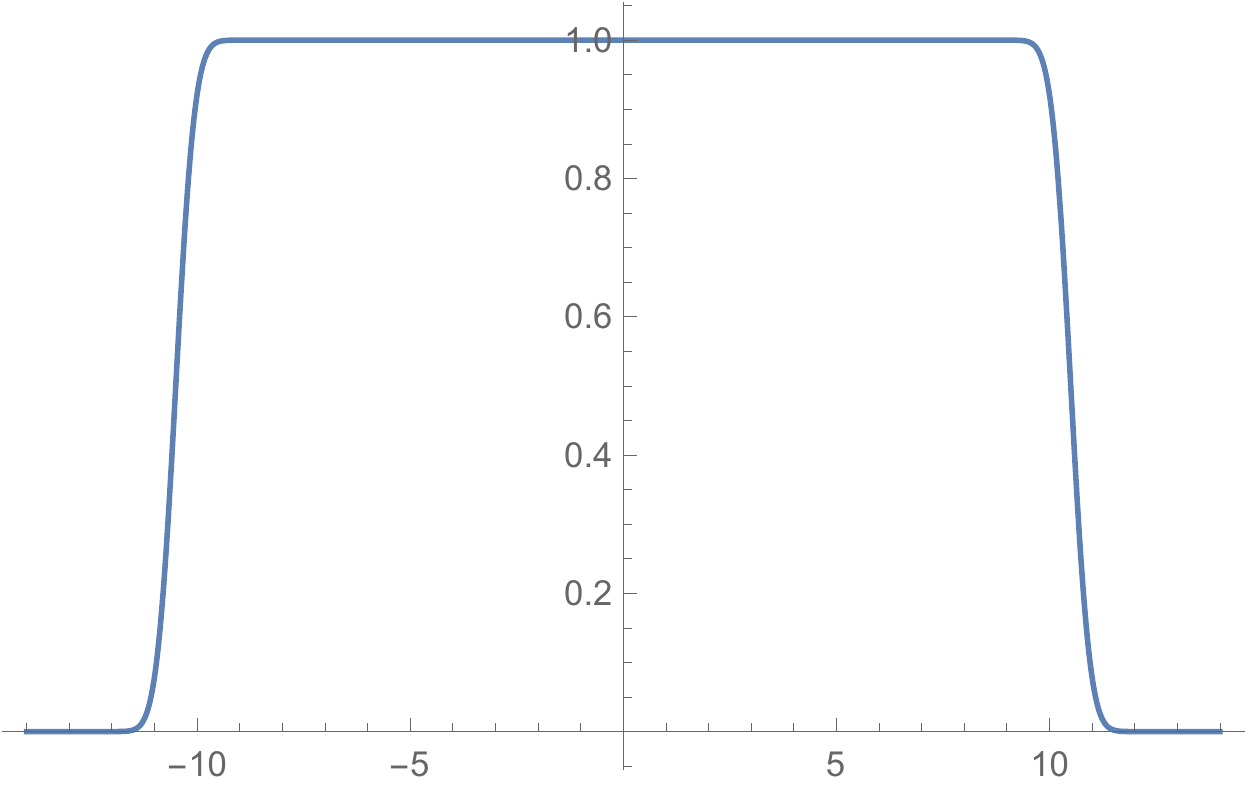}
\caption{Filter function in $-14<x<14$.
\label{fig:filter}}
\end{figure}

Notice that in Eq.~(\ref{eq:L2L}), because of the rapid decay of the weight function $e^{-t^2}$, 
one only needs to accurately approximate the function $h(t_l, t_r)$ in the region $[-7,7]^2$.
In our algorithm implementation, we define a filter function 
$$ \text{filter}(x,\epsilon)=\frac12 \left(\text{erf}(\frac{x/7+1.5}{\epsilon}) - \text{erf}(\frac{x/7-1.5}{\epsilon})\right) $$
where we set $\epsilon=\frac{1}{14}$ so that the function is approximately $\text{filter} \approx 1$ 
when $-7<x<7$ ($1-\text{filter}(7,\frac{1}{14}) = 2.09\text{e-}23$), and smoothly decays to $\text{filter} \approx 0$ at 
$\pm 14$ ($\text{filter}(14,\frac{1}{14}) = 2.09\text{e-}23$) , as shown in Fig.~\ref{fig:filter}.
At a leaf node, the integral is computed analytically either using 
$$\int_{a}^{b} e^{-x^2 -2 i x t} d x = \frac12 \sqrt{\pi} \left(
  \text{Fadd}(ia-t) e^{-a^2-2iat} - \text{Fadd}(ib-t) e^{-b^2-2ibt} \right) $$
for a boundary node, or 
$$\int_{a}^{b} e^{-2 i x t} d x = \frac{i}{2 t}  \left(
  e^{-2ibt} - e^{-2iat} \right) $$
for an interior node and then evaluated at a set of uniformly distributed $(2M)^2$
sample points in $[-14,14]^2$ for the two $t$-variables. The function values are then
filtered by the pointwise multiplication with the filter function for each variable. 
The Fourier series of the leaf node function, when needed, can be derived by a 
$2D$ FFT using the filtered function values. In the formula,  we use the Faddeeva function 
\cite{abrarov2012fourier,abrarov2011efficient,gautschi1970efficient,karbach2014decay} defined as
$\text{Fadd}(z)=e^{-z^2} \text{erfc}(-iz)$ for a complex number $z$, to avoid the possible overflow/underflow 
when computing small $e^{-t^2}$ times large $\text{erf}(a+it)$ values. 
An upward pass is then performed to recursively compute the parent's filtered function $h_p$ values
at the Fourier interpolation points using its children's filtered function values at different $t_l$, 
$t_m$, and $t_r$ interpolation points through 5 steps: (\rom{1}) multiplying two children's values 
at each sample point; (\rom{2}) point-wise multiplication with the filter function; (\rom{3}) applying the 
$1D$ fast Fourier transform (FFT) to the $t_m$ variable in the region 
$[-14,14]$ to get the $2M$ Fourier coefficients from the filtered function values at each $t_l,t_r$ interpolation
point; (\rom{4}) the parent's $h$-function value at each $t_l$ and $t_r$ interpolation point is derived 
by applying the formula 
$$ \frac{1}{\sqrt{\pi}} \int_{-\infty}^{\infty} e^{-t^2} e^{ i k \pi t /L} \text{d}t = e^{-\frac{k^2 \pi^2}{4 L^2}} $$ 
to integrate the Fourier series expansion of $t_m$ variable from (\rom{3}) analytically; 
and (\rom{5}) the function values will be further filtered. If needed, a $2D$ FFT can be performed to derive 
the parent's Fourier series expansion coefficients. Note that the Fourier series in the region 
$[-14,14]^2$ can be extended to the whole space $(-\infty, \infty)^2$ as such extension will only 
introduce an error within machine precision when evaluating the integral in Eq.~(\ref{eq:L2L}). 
At the root node, its $h$-function returns the $\phi$ value we are searching for.

The algorithm for efficiently evaluating Eq.~(\ref{eq:tridiag}) can be summarized as 
the following two passes. In the {\it downward pass}, the parent problem is decoupled by 
applying the Fourier transform to the coupling term, to obtain two child problems. 
At the finest level, a function with two $t$-variables is created for each leaf node 
followed by an {\it upward pass} to obtains each parent's function values at the Fourier 
interpolation points from those of its two children's functions. 
At the root level, the constant function (with null $t$-variables) gives the result of 
the integral in Eq.~(\ref{eq:tridiag}). 
The recursively implemented Matlab code for the {\it upward pass} is presented in Algorithm $1$.

\begin{table}[htbp]
\begin{mdframed}[
    linecolor=black,
    linewidth=2pt,
    roundcorner=4pt,
    backgroundcolor=olive!15,
    userdefinedwidth=\textwidth,
    nobreak=true,
]
{\bf function compute\_tri(inode)} 

\quad global NODES \%NODES contains the node informations. 

\quad if NODES(5,inode) == 0,  \% inode is a leaf node. 

\qquad leafnode(inode);   

\quad else 

\qquad    child1=NODES(5,inode); child2=NODES(6,inode); \% find children

\qquad    compute\_tri(child1); \% find child1's coefficients.

\qquad    compute\_tri(child2); \% find child2's coefficients.
   
\quad 

\qquad    \% combine children's coefficients to get parent's coefficients.

\qquad     if NODES(3,child1)==1 \&\& NODES(3,child2)==2

\qquad \quad        root(child1,child2); \% parent is the root node.

\qquad    elseif NODES(3,child1)==1 \&\& NODES(3,child2)==4

\qquad \quad         leftbdry(child1,child2); \% parent is a left boundary node.

\qquad    elseif NODES(3,child1)==3 \&\& NODES(3,child2)==2

\qquad \quad         rightbdry(child1,child2); \% parent is a right boundary node.

\qquad   else 

\qquad \quad         interiornode(child1,child2); \% parent is an interior node.

\qquad end

\quad end

return

end
\end{mdframed}
\caption{Algorithm 1:  Recursive Matlab function for evaluating Eq.~(\ref{eq:tridiag}): upward pass}
\end{table}

\subsection{Preliminary Numerical Results}
We present some preliminary results to demonstrate the accuracy and efficiency of the numerical algorithm
for the tridiagonal system in Eq.~(\ref{eq:tridiag}). In the numerical experiments, we set all
$a_i's$ to $-1$, $b_1=0.5$, $b_2=2$, and all other $b_i's$ to $+1$.  
We first study the accuracy of the algorithm. For $N=4$, we compute a reference solution
using Mathematica with $PrecisionGoal \rightarrow 30$ and $WorkingPrecision \rightarrow 60$, the result is 
$\phi= 2.2893342150887782603$.
For $N=8$, Mathematica returns the result $\phi=6.6242487478171897$ with an estimated
error $4.25\text{e-}5$, even though $PrecisionGoal \rightarrow 20$ and $WorkingPrecision \rightarrow 40$ are requested.
For $N>8$, direct computation using Mathematica simply becomes impossible.
In Table \ref{tab:triAccuracy}, we show the Matlab results for different dimensions $N$
and numbers of terms $2M$ in the Fourier series expansion. For all cases, our results
converge when $M$ increases. For $N=4$, our result matches Mathematica result to machine
precision as soon as enough Fourier terms are used. For $N=8$, our converged results agree 
with Mathematica result in the first $10$ digits, and we strongly believe our results are more accurate.
The numerical tests are performed on a laptop computer with Intel i7-3520M CPU @2.90GHz, with
$8.00G$ RAM.  For $N=1024$ and $M=512$, approximately $1024 \times 1024 \times 2047$ function
values at the Fourier interpolation points have to be stored in the memory ($\approx 16G$), 
which exceeds the installed RAM size, hence no result is reported.  
\begin{table}[htbp]
\centering
\begin{tabular}{|c|c|c|c|}
\hline
$N$ & 4 & 8 & 16  \\
\hline
$M$=16  & 2.326607912389402 & 6.736597967982384  & 56.44481808043047   \\
$M$=32  & 2.289334215119377 & 6.624246691958165  & 55.44625398858155   \\
$M$=64  & 2.289334215088778 & 6.624246691490006  & 55.44625397830180   \\
$M$=128 & 2.289334215088779 & 6.624246691490009  & 55.44625397830178  \\
$M$=256 & 2.289334215088778 & 6.624246691490005  & 55.44625397830176  \\
$M$=512 & 2.289334215088778 & 6.624246691490003  & 55.44625397830172  \\
\hline
\hline
$N$ & 32 & 64 & 1024  \\
\hline
$M$=16  & 3962.697712673563 & 19531008.87334120  & 1.182324449792241e+118   \\
$M$=32  & 3884.575992952042 & 19067179.07844248  & 1.019931849681238e+118   \\
$M$=64  & 3884.575991340509 & 19067179.06178229  & 1.019931834748418e+118   \\
$M$=128 & 3884.575991340506 & 19067179.06178229  & 1.019931834748411e+118   \\
$M$=256 & 3884.575991340500 & 19067179.06178224  & 1.019931834748369e+118   \\
$M$=512 & 3884.575991340498 & 19067179.06178220  & N/A   \\
\hline
\end{tabular}
\caption{Computed $\phi$ values for different dimensions and number of Fourier terms.
\label{tab:triAccuracy}}
\centering
\end{table}

We demonstrate the efficiency of our algorithm by presenting the Matlab simulation
time for different dimensions. In the experiment, we present the CPU times for different
$M$ and $N$ values, and the unit is in seconds. Clearly, the CPU time grows approximately 
linearly as the dimension $N$ increases. As a 3-variable \{$t_l$, $t_m$, and $t_r$\} 
function has to be processed in the current implementation when finding the parent's 
values at the Fourier interpolation points, the CPU time grows approximately by a 
factor of $8$ as $M$ doubles. For $N=2048$ and $M=256$, approximately $8G$ memory is 
required, which exceeds the maximum available RAM size, hence no result is reported.
\begin{table}[htbp]
\centering
\begin{tabular}{|c|c|c|c|c|c|}
\hline
$N$ & 4 & 8 & 16 & 32 & 64 \\
\hline
$M=32$ CPU time & 0.02 & 0.04 & 0.15 & 0.34 & 0.80 \\ 
\hline
$M=64$ CPU time & 0.06 & 0.27 & 0.96 & 2.26 & 4.60 \\ 
\hline
$M=128$ CPU time & 0.12 & 2.15 & 6.81 & 16.5 & 39.1 \\ 
\hline
$M=256$ CPU time & 0.24 & 12.5 & 49.1 & 130 &  315 \\ 
\hline
\hline
$N$ & 128 & 256 & 512 & 1024 & 2048 \\
\hline
$M=32$ CPU time & 1.95 & 3.75 & 8.29 & 16.6 & 40.1 \\ 
\hline
$M=64$ CPU time & 10.2 & 23.5 & 49.2 & 103 & 214 \\ 
\hline
$M=128$ CPU time & 88.5 & 228 & 445 & 970 &  1942 \\ 
\hline
$M=256$ CPU time & 718 & 1769 & 3476 & 8310 & N/A \\ 
\hline
\end{tabular}
\caption{CPU time (in seconds) for different $N$ and $M$ values.
\label{tab:triEfficiency}}
\centering
\end{table}

\section{Case II: Exponential Matrix}
\label{sec:rankone}
In the second case, we consider a matrix $A$ defined by the {\it exponential covariance function} 
\[ A_{i,j}=e^{-|z_i-z_j|/\beta}, \beta > 0.
\] 
To simplify the discussions, we consider a simple 1D setting from spatial or temporal statistics 
and assume that the rate of decay $\beta = 1$ and each random number $X_j$ is observed at a location 
$z_j \in [0,b_z]$. We assume the $z$-locations $\{z_j \in [0,b_z]$, $j=1, \ldots, N=2^P \}$ are 
ordered from smallest to largest and the matrix entries are ordered accordingly. 
We demonstrate how to evaluate the $N$-dimensional integral
\begin{equation}
\phi(\ba, \bb; A)= \int_{\ba}^{\bb} f(\bx | A) \text{d}\bx= 
  \int_{a_1}^{b_1} \cdots \int_{a_N}^{b_N} 
  \exp \left(-\frac12 \bx^T A \bx \right) \text{d}x_N \cdots \text{d}x_1,
\end{equation}
for the given constant vectors $\ba$ and $\bb$ using $O(N)$ operations. Results for 
different $\beta$ values can be derived by rescaling the $z$-locations and $x$-variables. 
The presented algorithm can be easily generalized to 
$\int_{\ba}^{\bb} H(\bx) f(\bx | A) \text{d}\bx$ when $H(\bx)$ is a low-rank function.
 
Similar to the tridiagonal matrix case, we generate a binary tree by recursively dividing the 
parent's $z$-location set (or equivalently the $x$-variable set) into two child subsets, 
each containing exactly 
half of its parent's points. The hierarchical binary tree is then reflected as a hierarchical matrix 
as demonstrated in Fig.~\ref{fig:hmatrix}. Unlike the (uniform) binary tree generated for the 
$z$-location set, the corresponding structure in the matrix sub-division process can be considered as an adaptive quad-tree, 
where only the diagonal blocks of the matrix are subdivided. Once an off-diagonal block is generated, 
it becomes a leaf node and no further division is required. Because of the hierarchical structure of
the matrix and the low-rank properties of the off-diagonal blocks (which will be discussed next), 
the exponential matrix is a special $\mathcal{H}$-matrix.

\subsection{Divide and Conquer on a Hierarchical Tree}
Unlike the tridiagonal system, each off-diagonal matrix in this case is a dense matrix. For this
exponential matrix, all the off-diagonal matrices are rank-1 matrices, which can be seen from the 
separation of variables
\[ e^{-|z-y|} = \left\{ \begin{array}{l} e^{-z} e^y, \quad z \geq y \\ e^z e^{-y}, \quad z < y \end{array}
  \right. 
\]
In matrix language, the off-diagonal block $A_{1,3}$ can be written as 
\begin{equation}
\label{eq:rank1} 
  [A_{1,3}(y_i, z_j)] = [e^{-y_{N/2+1}}, \ldots, e^{-y_{N}}]^T [e^{z_1},\ldots, e^{z_{N/2}}] 
\end{equation}
for $i=N/2+1, \ldots, N$ and $j=1, \ldots, N/2$.
The singular value decomposition of $A_{1,3}$ can be easily derived using Eq.~(\ref{eq:rank1}) as
\begin{displaymath}
A_{1,3}= \bu \lambda \bv^T
\end{displaymath}
where the left and right singular vectors $\bu$ and $\bv$ are of size $\frac{N}{2} \times 1$ and 
respectively the normalized vectors of the discretized functions $e^{-y}$ and $e^z$. 

When the $x$-variables are divided into $2$ subsets $\bx_{1,1}$ and $\bx_{1,2}$, the root matrix $A$ 
can be subdivided accordingly into $4$ blocks
\begin{equation*}
A = 
\begin{bmatrix}
A_{1,1} & A_{1,2}=\bv \lambda \bu^T \\
A_{1,3}=\bu \lambda \bv^T & A_{1,4}
\end{bmatrix}, 
\bx = 
\begin{bmatrix}
\bx_{1,1}\\
\bx_{1,2}
\end{bmatrix}, 
\end{equation*}
where the first index of $A_{i,j}$ is the current level of the block matrix and the 
second index is its order in this level. Same indexing rules are used for the $z$-locations 
and $x$-variables. Completing the square, the quadratic form in the integrand can be reformulated as
\begin{displaymath}
\begin{array}{rcl}
\bx^T A \bx &=& \bx_{1,1}^T A_{1,1} \bx_{1,1} + {\color{blue} \bx_{1,1}^T \bv \lambda \bu^T \bx_{1,2} 
  + \bx_{1,2}^T \bu \lambda \bv^T \bx_{1,1} } + \bx_{1,2}^T A_{1,4} \bx_{1,2} \\
                &=& \bx_{1,1}^T A_{1,1} \bx_{1,1} + \bx_{1,2}^T A_{1,4} \bx_{1,2} + 
			{\color{red} \left( (\gamma \bu^T \bx_{1,2} + \frac{1}{\gamma} \bv^T \bx_{1,1}) \sqrt{\lambda} \right)^2 
                     }   \\
                & &   -{\color{blue} \bx_{1,2}^T \gamma^2 \bu \lambda \bu^T\bx_{1,2} } - 
                         {\color{blue} \bx_{1,1}^T \frac{1}{\gamma^2} \bv \lambda \bv^T \bx_{1,1} } \\
                &=& {\color{green} \bx_{1,1}^T (A_{1,1} - \frac{1}{\gamma^2} \bv \lambda \bv^T) \bx_{1,1}} + 
                    {\color{green} \bx_{1,2}^T (A_{1,4} - \gamma^2 \bu \lambda \bu^T) \bx_{1,2} }  \\
                & & 	+ {\color{red} \left( (\gamma \bu^T \bx_{1,2} + \frac{1}{\gamma} \bv^T \bx_{1,1}) \sqrt{\lambda} \right)^2  } \\
\end{array}
\end{displaymath}
where the first two {\color{green} green} terms are the child problems to be processed recursively at finer levels in
the divide-and-conquer strategy, $\gamma$ is a constant to be determined, and the last {\color{red} red} term 
shows how the two child problems are coupled. Similar to the tridiagonal case,
by introducing a single $t$-variable and applying the Fourier transform formula in Eq.~(\ref{eq:Fourier})  
to the coupling term (in red), we get
\begin{displaymath}
\begin{array}{l}
 \int_{\ba}^{\bb} 
  e^{ - \frac12 \bx^T A \bx} \text{d}\bx = c \int_{-\infty}^{\infty} e^{- t^2} h_{1,1}(t) h_{1,2}(t) \text{d}t 
\end{array}
\end{displaymath}
where $h_{1,1}(t)$ and $h_{1,2}(t)$ are the single $t$-variable functions for the two child nodes given by
\begin{equation}
	\label{eq:rootdivide}
\begin{array}{l}
 h_{1,1}(t) = \int_{a_1}^{b_1} \cdots \int_{a_k}^{b_k}  
  e^{-  \frac12 \bx_{1,1}^T (A_{1,1} - \frac{\lambda}{\gamma^2} \bv \bv^T) \bx_{1,1} + 
     i t   \frac{\sqrt{2 \lambda}}{\gamma} \bv^T \bx_{1,1} } \text{d}\bx_{1,1},\\
 h_{1,2}(t) = \int_{a_{k+1}}^{b_{k+1}} \cdots \int_{a_N}^{b_N}  
  e^{ -  \frac12 \bx_{1,2}^T (A_{1,4} -\gamma^2 \lambda  \bu \bu^T) \bx_{1,2} -    
    i t \sqrt{2 \lambda} \gamma \bu^T \bx_{1,2}  } \text{d}\bx_{1,2}.
\end{array}
\end{equation}
Note that the $x$-variables are completely decoupled in the two child problems, and the coupling is now 
through the $t$-variable.

In order to have a divide-and-conquer algorithm on the hierarchical tree structure, the two child problems
should have the following properties:
\begin{itemize}
	\item By properly choosing the parameter $\gamma$, the new matrices 
		$A_{1,1} - \frac{\lambda}{\gamma^2} \bv \bv^T$ and $A_{1,4} -\gamma^2 \lambda  \bu \bu^T$
		should be symmetric positive definite; and 
	\item The off-diagonal blocks of these new matrices should be low-rank. 
\end{itemize}
We found that the choice of $\gamma$ is not unique, and there exist a range of $\gamma$ values for the child
problems to have these properties. The choice of $\gamma$ is addressed next.

\subsection{Potential Theory based Analysis}
In this section, we apply the potential theory from the analysis of ordinary and partial differential 
equations and show how the divide-and-conquer strategy can be successfully performed on the 
hierarchical tree structure. Purely numerical linear algebra based approaches for more general cases will 
be briefly addressed later.

\subsubsection{Green's Functions}
We present the results for $b_z=1$ to simplify the notations and assume $z_j \in [0,1]$. We start from the 
observation that 
$$G(z,y)=\frac12 e^{-|z-y|} = \left\{ 
\begin{array}{l} 
	coef \cdot g_r(z) \cdot g_l(y), \quad z\geq y, \\
	coef \cdot g_r(y) \cdot g_l(z), \quad z<y
\end{array}
\right.
$$ is the domain Green's function of the ordinary differential
equation (ODE) two-point boundary value problem
\begin{equation}
	\label{eq:rootode}
	\left\{
	\begin{array}{l}
		u(z)-u''(z)=f(z), \quad z \in [0,1], \\
		u(0)=u'(0), \quad u(1)=-u'(1),
	\end{array}
	\right.
\end{equation}
where $coef=\frac12$, $g_l(z) = e^{z-1}$ and $g_r(z)=e^{1-z}$.
The proof is simply a straightforward validation that $u(z)= \int_0^1 G(z,y) f(y) dy$ satisfies both the
ODE and boundary conditions. 

In the following discussions, we consider the continuous version of the
original matrix problem, where the matrix $A$ is the discretized Green's function $G(z,y)$,
the two off-diagonal submatrices $A_{1,2}$ and $A_{1,3}$ are the discretized 
$g_r(z) \cdot g_l(y)$ and $g_r(y) \cdot g_l(z)$, respectively. 
Some simple algebra manipulations show that the submatrices $A_{1,1}- \frac{\lambda}{\gamma^2} \bv \bv^T$ 
and $A_{1,4}-\gamma^2 \lambda  \bu \bu^T$ can be considered as the discretized 
$G(z,y)-\tilde{\gamma}^2 \cdot g_l(z) \cdot g_l(y)$ and 
$G(z,y)-\frac{1}{\tilde{\gamma}^2} \cdot g_r(z) \cdot g_r(y)$, and the coefficients 
$i t   \frac{\sqrt{2\lambda}}{\gamma} \bv^T$  and $i t \sqrt{2\lambda} \gamma \bu^T$  for the linear 
terms of the $x$-variables $\bx_{1,1}$ and $\bx_{1,2}$ in Eq.~(\ref{eq:rootdivide}) are the discretized 
$i t \tilde{\gamma} g_l(z)$ and $ i t \frac{1}{\tilde{\gamma}} g_r(z)$, respectively.

\vspace{0.1in}
{\noindent \bf Remark:} The observation also allows easy proof of the positive definiteness of the
matrix $A$, which is the discretized Green's function $G(z,y)$. 
In order to show that for any vector $\bff \neq \bf{0}$, the quadratic form satisfies $\frac12 \bff^T A \bff >0$, 
we consider its continuous version defined as 
$$\int_0^1 f(z) \left(\int_0^1 G(z,y) f(y) dy \right) dz = \int_0^1 f(z) u(z) dz$$ 
where $u(z)= \int_0^1 G(z,y) f(y) dy$ and $f(y)$ is the continuous version of the (discretized) 
vector $\bff$. As $f(z)=u(z)-u''(z)$, applying the integration by parts, we have
$$ \int_0^1 f(z) u(z) dz = \int_0^1 \left( u^2(z) + (u'(z))^2\right) dz 
  -u'(1) \cdot u(1) + u'(0) \cdot u(0). $$ 
As $f(z) \neq 0$, therefore $u(z) \neq 0$, and applying the boundary conditions of the ODE, we have 
$\int_0^1 f(z) u(z) dz>0$. We refer to the two-variable function $G(z,y)$ as a positive definite 
function.  The positive definiteness of the matrix $A$ can be proved in a similar way using 
the {\it discretized} integration by parts.

A particular choice of $\gamma$ can be determined by considering the corresponding 
child ODE problems as follows. We first study the root problem and define its two children as 
the {\it left child} and {\it right child}, and the locations $z_i$ of the left child and $z_j$ 
of the right child satisfy the condition $z_i< z_j$ as the $z$-locations of the $x$-variables 
in the two child problems are separated and ordered. We pick a location $\zeta$ between the two
clusters of $z$-locations. Note that the choice of $\zeta$ is not unique. We have the following
results for the root node.
\begin{theorem}
If we choose $\tilde{\gamma} = \frac{e^{-\zeta}}{\sqrt{2}}$, then 
\begin{itemize}
	\item for the left child, the new function $G_l(z,y)=G(z,y)-\tilde{\gamma}^2 \cdot g_l(z) \cdot g_l(y)$
	  is the Green's function of the ODE problem
	\begin{displaymath}
	\left\{
	\begin{array}{l}
		\bu_1(z)-\bu_1''(z)=f(z), \quad z \in [0,\zeta], \\
		\bu_1(0)=\bu_1'(0), \quad \bu_1(\zeta)=0.
	\end{array}
	\right.
	\end{displaymath}
		The function $G_l(z,y)$ is {\it positive definite}.
	\item For the right child, the new function 
		 $G_r(z,y)=G(z,y)-\frac{1}{\tilde{\gamma}^2} \cdot g_r(z) \cdot g_r(y)$
  	is the Green's function of the ODE problem
	\begin{displaymath}
	\left\{
	\begin{array}{l}
		\bu_2(z)-\bu_2''(z)=f(z), \quad z \in [\zeta, 1], \\
		\bu_2(\zeta)=0, \quad \bu_2(1)=-\bu_2'(1).
	\end{array}
	\right.
	\end{displaymath}
		The function $G_r(z,y)$ is {\it positive definite}.
	\item The two child ODE problem solutions $\bu_1(z)$ and $\bu_2(z)$ can be derived by subtracting a
		single layer potential defined at $z=\zeta$ from the parent's solution $u(z)$ of 
		Eq.~(\ref{eq:rootode}), so that solutions $\bu_1(z)$ and $\bu_2(z)$ satisfy the zero 
		interface condition at $z=\zeta$. The other boundary condition for each child ODE 
		problem is the same as its parent's boundary condition.
\end{itemize}
\end{theorem}
These results can be easily validated by plugging in the functions to the ODE problems. The positive definiteness
of the child Green's function can be proved using the same integration by part technique as we did for the
parent's Green's function.

For a general parent node on the tree structure, we have the following generalized results.
\begin{theorem}
	Consider a parent node with the corresponding function $G_p(z,y)$ defined on the interval
        $[a,b]$, and $\zeta$ is a point separating the two children's $z$-locations. Then there
	exists a number $\tilde{\gamma}$ which depends on $\zeta$, such that 
\begin{itemize}
	\item for the left child, the new function $G_l(z,y)=G(z,y)-\tilde{\gamma}^2 \cdot g_l(z) \cdot g_l(y)$
	  is the Green's function of the ODE problem
	\begin{displaymath}
	\left\{
	\begin{array}{l}
		\bu_1(z)-\bu_1''(z)=f(z), \quad z \in [a,\zeta], \\
		\mbox{same boundary condition as parent at $x=a$, and } \bu_1(\zeta)=0.
	\end{array}
	\right.
	\end{displaymath}
		The function $G_l(z,y)$ is {\it positive definite}.
	\item For the right child, the new function 
		 $G_r(z,y)=G(z,y)-\frac{1}{\tilde{\gamma}^2} \cdot g_r(z) \cdot g_r(y)$
  	is the Green's function of the ODE problem
	\begin{displaymath}
	\left\{
	\begin{array}{l}
		\bu_2(z)-\bu_2''(z)=f(z), \quad z \in [\zeta, b], \\
		\bu_2(\zeta)=0, \mbox{ and same boundary condition as parent at $x=b$}.
	\end{array}
	\right.
	\end{displaymath}
		The function $G_r(z,y)$ is {\it positive definite}.
	\item The two child ODE problem solutions $\bu_1(z)$ and $\bu_2(z)$ can be derived by subtracting a
		single layer potential defined at $z=\zeta$ from the parent's solution $u(z)$ of 
		Eq.~(\ref{eq:rootode}), so that solutions $\bu_1(z)$ and $\bu_2(z)$ satisfy the zero 
		interface condition at $z=\zeta$. The other boundary condition for each child ODE 
		problem is the same as its parent's boundary condition.
\end{itemize}
\end{theorem}
The detailed formulas for the number $\tilde{\gamma}$ and Green's functions are presented in the Appendix. 
The proof of the theorem is simply validations of the formulas.

\subsubsection{Parent-children Relations} 
In the matrix form,
for a general parent node at level $l$ in the hierarchical tree structure with left child $1$ and 
right child $2$, its $h$-function
\begin{equation}
\label{eq:gfun}
	h_p (\bt_{p}) = \int_{\ba_p}^{\bb_p} e^{ - \frac12 \bx_p^T A_p \bx_p} e^{ i \bt_{p}^T D_p \bx_p } \text{d}\bx_p 
\end{equation}
can be decomposed into two child problems as 
\[
	h_p(\bt_{p}) = \frac{1}{\sqrt{\pi}} \int_{-\infty}^{\infty} e^{-t_{new}^2 } h_{1}(\bt_1) h_{2}(\bt_2) \text{d}t_{new},
\]
where 
\begin{equation}
\label{eq:decoup2}
\begin{array}{l}
	h_{1}(\bt_1) = \int_{\ba_1}^{\bb_1} e^{-\frac12 \bx_1^T A_1 \bx_1} e^{i t_{new} \tilde{\gamma} \bg_l^T(\bz_1) \cdot \bx_1} 
	   e^{ i \bt_p^T D_{p,1} \bx_1 } \text{d}\bx_1 = 
	\int_{\ba_1}^{\bb_1} e^{-\frac12 \bx_1^T A_1 \bx_1} e^{ i \bt_1^T D_{1} \bx_1 } \text{d}\bx_1 , \\
	h_{2}(\bt_2) = \int_{\ba_2}^{\bb_2} e^{-\frac12 \bx_1^T A_2 \bx_1} e^{i t_{new} \frac{1}{\tilde{\gamma}} \bg_r^T(\bz_2) \cdot \bx_2} 
	  e^{ i \bt_p^T D_{p,2} \bx_2 } \text{d}\bx_2 = \int_{\ba_2}^{\bb_2} e^{-\frac12 \bx_2^T A_2 \bx_2} e^{ i \bt_2^T D_{2} \bx_2 } \text{d}\bx_2.
\end{array}
\end{equation}
In the formulas, $\bt_p$ is the vector containing all the $t$-variables introduced at coarser levels to subdivide $p$'s 
parents' $h$-functions. $\bx_p= [ \bx_1 ; \bx_2]$, $\bx_1$, and $\bx_2$ 
are respectively the vectors containing the $x$-variables of the parent $p$, child $1$, and child $2$.  
\{$\ba_p$,$\bb_p$\}, \{$\ba_1$,$\bb_1$\}, and \{$\ba_2$,$\bb_2$\} are respectively the lower and 
upper integration bounds of $\bx_p$, $\bx_1$, and $\bx_2$. $A_p= \left[ \begin{array}{cc} A_{l,1} & A_{l,2} \\ A_{l,3} & A_{l,4} \end{array} \right]$, $A_1$, and $A_2$ are
respectively the discretized Green's functions of the parent $p$, child $1$, and child $2$, which satisfy
$$ A_1=A_{l,1}- \tilde{\gamma}^2 \cdot \bg_l(\bz_1) \cdot \bg_l^T(\bz_1), \quad 
A_2=A_{l,4}- \frac{1}{\tilde{\gamma}^2} \cdot \bg_r(\bz_2) \cdot \bg_r^T(\bz_2),$$ 
$\bz_p=[\bz_1; \bz_2]$, $\bz_1$, and $\bz_2$ are respectively the $z$-location vectors of the parent $p$ and 
child $1$ and $2$, and $\bg_l(\bz_1)$ and $\bg_r(\bz_2)$ are the discrete function values of $g_l(z)$ and $g_r(z)$ 
in the Green's functions evaluated at different $z$-locations.
$\bt_{p}^T D_p \bx_p$ is a scalar term representing the linear combinations of the $t_k \cdot x_j$ terms, and 
by separating the $x$-variables, it can be written as 
$$ \bt_{p}^T D_p \bx_p = \bt_{p}^T D_{p,1} \bx_1 + \bt_{p}^T D_{p,2} \bx_2.$$
After introducing the new $t$-variable $t_{new}$ to divide the parent's problem to two subproblems of child $1$ 
and child $2$, each with half of the parent $p$'s $x$-variables, we have 
$\bt_1=[\bt_p; t_{new}]$, $\bt_2=[\bt_p; t_{new}]$, and
\begin{equation} 
	\label{eq:parentchild01}
	\left\{ 
	\begin{array}{l} \bt_1^T D_{1} \bx_1 = \bt_{p}^T D_{p,1} \bx_1 +  t_{new} \tilde{\gamma} \bg_l^T(\bz_1) \cdot \bx_1, \\
 \bt_2^T D_{2} \bx_2  =  \bt_{p}^T D_{p,2} \bx_2 +  t_{new} \frac{1}{\tilde{\gamma}} \bg_r^T(\bz_2) \cdot \bx_2.
	\end{array}
	\right.
\end{equation}

For the root node, $A_p$ is the given matrix $A$ and $\bt_p$ is an empty set. At a leaf node, we have
$$h_{leaf}(\bt_{leaf}) = \int_{a_k}^{b_k} e^{-\frac12 \alpha_k x_k^2} e^{i (\bt_{leaf}^T D_{leaf}) x_k} \text{d}x_k $$
where $D_{leaf}$ is a column vector of the same size as $\bt_{leaf}$ (the size equals to the number of levels
in the hierarchical tree structure). Analytical formula is available for $h_{leaf}(\bt_{leaf})$ using 
\begin{equation}
\label{eq:leaf}
	\int_a^b e^{-x^2} e^{-2 i t x} \, \text{d} x = \frac{1}{2} \sqrt{\pi } e^{-t^2} (\text{erf}(b+i t)-\text{erf}(a+i t)).
\end{equation}

\subsubsection{Dimension Reduction and Effective Variables}
Note that for a node at level $l$, its $h$-function $h(\bt)$ will contain as many as $l$ $t$-variables introduced at 
parent levels. Therefore for
a $N$-dimensional problem, the number of $t$-variables for a leaf node can be as many as $\log (N)$. However, inspecting
the term $(\bt_{leaf}^T D_{leaf}) x_k$ for the function $h_{leaf} (\bt_{leaf})$, if one introduces a new single variable 
$w = \bt_{leaf}^T D_{leaf}$, then $h_{leaf}$ is effectively a {\it single} variable function of $w$. We therefore study the 
{\it effective} variables and their properties in this section.

From Eq.~(\ref{eq:parentchild01}), we see that when a new $t$-variable $t_{new}$ is introduced to divide the parent problem
into two child problems, the additional terms added to the linear terms of the $x$-variables in the exponent are
$t_{new} \tilde{\gamma} \bg_l^T(\bz_1) \cdot \bx_1$ for child $1$ and $t_{new} \frac{1}{\tilde{\gamma}} \bg_r^T(\bz_2) \cdot \bx_2$ 
for child $2$, where $ \bg_l(\bz_1)$ and $\bg_r(\bz_2)$ are the discrete function values of $g_l(z)$ and $g_r(z)$ 
in the Green's functions evaluated at different $z$-locations. For all the Green's functions, 
$g_l(z)$ and $g_r(z)$ are always a combination of the basis functions $e^z$ and $e^{-z}$. This can be seen either from
the ODE problems or from the Green's functions in the Appendix. Therefore, switching the basis to $e^z$ and $e^{-z}$,
the term $ \bt^T D \bx $ can always be written as 
\begin{equation}
	\bt^T D \bx = (w_1 e^{\bz} + w_2 e^{-\bz})^T \cdot \bx,
\end{equation}
where $e^{\bz}$ and $e^{-\bz}$ are the vectors derived by evaluating the functions $e^z$ and $e^{-z}$ at the
$z$-locations. Clearly, after this change of variables from $t$-variables to \{$w_1$, $w_2$\}, each $h$-function is
effectively a function with no more than $2$ variables. We define $w_1$ and $w_2$ as the effective $w$-variables. 

Our numerical experiments show that at finer levels of the hierarchical tree structure when the interval size
of the tree node becomes smaller, the two basis functions $e^z$ and $e^{-z}$ are closer to linear dependent 
which will cause numerical stability issues. For better stability properties, orthogonal or near orthogonal
basis functions are used. A sample basis is $\{ \Phi_1(z)=\cosh(z-c), \Phi_2(z)=\frac{\sinh(z-c)}{b-a} \}$ when
the $z$-locations of the $x$-variables are in the interval $[a,b]$. When $c$ is the center of the interval, 
the two functions are orthogonal to each other when measured using the standard $L_2$ norm with a constant weight function.
For a parent node with effective $w$-variables $\{w_1^p, w_2^p\}$ and basis functions $\{ \bPhi_1^p, \bPhi_2^p \}$,
where the vector $\bPhi$ represents the discretized $\Phi(z)$ at the $z$-locations, in the divide-and-conquer strategy, 
the effective $w$-variables should satisfy the relations 
\begin{equation}
\label{eq:pc-effective}
  \left\{
	\begin{array}{l}
	w_1^p \Phi_1^p + w_2^p \Phi_2^p +  t_{new} \tilde{\gamma} g_l(z) = w_1^1 \Phi_1^1 + w_2^1 \Phi_2^1, \\
	w_1^p \Phi_1^p + w_2^p \Phi_2^p +  t_{new} \frac{1}{\tilde{\gamma}} g_r(z) = w_1^2 \Phi_1^2 + w_2^2 \Phi_2^2,
	\end{array}
\right.
\end{equation}
where $\{\Phi_1^p, \Phi_2^p\}$, $\{\Phi_1^1,\Phi_2^1\}$, $\{\Phi_1^2,\Phi_2^2\}$ are 
respectively the continuous basis for the parent, child 1, and child 2, 
$\{w_1^1, w_2^1\}$ and $\{w_1^2, w_2^2\}$ are the effective $w$-variables of child $1$ 
and child $2$ for the discrete basis vectors $\{ \bPhi_1^1, \bPhi_2^1 \}$ and $\{ \bPhi_1^2, \bPhi_2^2 \}$, respectively. 
In the Appendix, we
present the detailed formulas demonstrating the relations between parent $p$'s and children's effective $w$-variables
for the basis choice \{$\cosh(z-c)$, $\frac{\sinh(z-c)}{b-a}$\}.

In the tridiagonal case discussed in Section \ref{sec:tridiag}, we only need to study the $h$-functions
when their $t$-variables satisfy $|t_j|<7$, as outside the interval the integrand value is controlled by the  
factor $e^{-t_j^2}$ and hence can be neglected. Similar results can be obtained for the exponential case,
when a proper set of basis is chosen. Assuming all the $z$-locations are approximately uniformly distributed in
the interval $[0,1]$, we have the following theorem for the effective $w$-variables $w_1$ and $w_2$.
\begin{theorem}
Assume the $N \times N$ matrix $A$ is defined by the exponential covariance function, the $z$-locations are 
uniformly distributed in the interval $[0,1]$, and all the $t$-variables satisfy $|t_j|<7$. When the basis
functions are chosen as $\{ \Phi_1(z)=\cosh(z-c), \Phi_2(z)=\frac{\sinh(z-c)}{b-a} \}$ for each tree node, then
there exists a constant $C$ independent of $N$, such that the corresponding effective $w$-variables $w_1$ and 
$w_2$ (combinations of the $t$-variables) satisfy the conditions $|w_1| \leq C$ and $|w_2| \leq C$.
\end{theorem}

The proof of this theorem is simply the leading order analysis of the parent-children effective $w$-variable 
relations, and the fact that $\cos(h) = 1 + \frac{h^2}{2} + O\left(h^4\right)$, 
$\frac{\sinh(h)}{2h} = \frac{1}{2}+\frac{h^2}{12}+O\left(h^4\right)$, 
$\sinh \left(\frac{h}{2}\right) \sqrt{\sinh (h) \text{csch}(2 h) \text{csch}(h)}=
  \frac{\sqrt{h}}{2 \sqrt{2}}-\frac{7 h^{5/2}}{48 \sqrt{2}}+O\left(h^{9/2}\right)$, 
and $\sum_{k=0}^L \sqrt{\frac{1}{2^k}} < \sqrt{2}+2$, where $L$ is the number of levels in the tree structure. 
We skip the proof details. Interested readers can request a copy of our Mathematica file for further 
details. We point out that when the basis functions are chosen as $\{ e^z, e^{-z}\}$, the effective $w$-variables
become unbounded.

\vspace{0.1in}
{\noindent \bf Remark:} In the numerical implementation, instead of using the upper bound $C$ for a tree node $j$, 
the ranges $C_1^j$ and $C_2^j$  of the effective $w$-variables $w_1$ and $w_2$ are computed using the 
parent-children effective $w$-variable relations in Eq.~(\ref{eq:pc-effective}) and stored in the memory. 
Similar to the tridiagonal case, a filter function is applied to the $h$-functions so that the filtered 
function smoothly decays to zero in the region $|w_1| \in [C_1, 2C_1]$ or $|w_2| \in [C_2, 2C_2]$, 
see Fig.~\ref{fig:filter}. Then the Fourier series of the filtered $h$-function is constructed in the 
region $[-2C_1, 2C_1] \times [-2 C_2, 2 C_2]$, and finally the constructed Fourier series is expanded to the
whole space when deriving parent's $h$-function values. In the algorithm implementation, when the 
uniform FFT \cite{frigo1998fftw} can no longer be applied, we use the open source NUFFT package 
developed in \cite{greengard2004accelerating,lee2005type} to accelerate the computation of the 
Fourier series.

\subsection{Pseudo-algorithm}
Similar to the tridiagonal case, the algorithm can be summarized as the following two passes: 
In the {\it downward pass}, the parent problem is decoupled by applying the Fourier transform 
to the coupling term, to obtain two child problems. Six coefficients $\{ c_1, c_2, c_3, c_4, c_5, c_6 \}$ 
are derived so that the effective $w$-variables of the current node satisfy
\begin{equation}
\label{eq:pc6number}
w_1= c_1 w_1^p + c_2 w_2^p + c_3 t_{new}, \quad w_2= c_4 w_1^p + c_5 w_2^p +c_6 t_{new},
\end{equation}
where $w_1^p$ and $w_2^p$ are the parent's effective $w$-variables. Also, the ranges $C_1$ and $C_2$ of 
the effective $w$-variables $w_1$ and $w_2$ are computed. A total of $8$ numbers are stored 
for each node. Note that both the storage and number of operations are constant for each tree
node. The pseudo-algorithm is presented in Algorithm $2$, where the details of computing the
$8$ numbers for each node is presented in the Appendix.

\begin{table}[htbp]
\begin{mdframed}[
    linecolor=black,
    linewidth=2pt,
    roundcorner=4pt,
    backgroundcolor=olive!15,
    userdefinedwidth=\textwidth,
    nobreak=true,
]
{\bf function compexp\_downward(inode)} 

\quad global NODES \%NODES contains the node informations. 

\quad global TRANSCoef \%TRANSCoef contains the $8$ numbers. 

\quad if NODES(5,inode) == 0,  \% inode is a leaf node. 

\qquad return;   

\quad else 

\qquad    child1=NODES(5,inode); child2=NODES(6,inode); \% find children
    
\qquad 	  compute the $8$ numbers using the formulas in Appendix for inode 

\qquad \qquad is a root, left boundary, right boundary, or interior node.

\qquad    compexp\_downward(child1); \% find child 1's $8$ numbers.

\qquad    compexp\_downward(child2); \% find child 2's $8$ numbers.
   
\quad end

return

end
\end{mdframed}
\caption{Algorithm 2: Recursive Matlab function for exponential case: downward pass}
\end{table}

At the finest level, a function with one effective variable is constructed analytically
using Eq.~(\ref{eq:leaf}). A numerically equivalent two-variable $\{w_1^{leaf}, w_2^{leaf} \}$
Fourier series expansion is then constructed by evaluating the analytical solution
at the interpolation points, applying the filter function, and then applying FFT to derive
the $2D$ Fourier series expansion which is considered valid in the whole space. 
An {\it upward pass} is then performed, to obtains each parent's Fourier coefficients from 
those of its two children's functions. For each parent node, we first replace the 
child's effective $w$-variables with $w_1^p$, $w_2^p$ and $t_{new}$ using 
Eq.~(\ref{eq:pc6number}) and the $6$ numbers from the downward pass, then evaluate 
each child's {\it global} Fourier series at the uniform interpolation points of
$w_1^p$, $w_2^p$ and $t_{new}$ (determined by the ranges $C_1$ and $C_2$ 
from the downward pass, we set the range of $t_{new}$ to $7$). In this step, we have
to use the NUFFT as the $8$ numbers for different tree nodes are different so 
the uniform FFT is not applicable. Multiplying the two children's function values and the 
filter function values at each interpolation point, we then apply the FFT to the 
$t_{new}$ variable and derive the Fourier series of $t_{new}$ at each $w_1^p$ and $w_2^p$ 
interpolation point. The integral 
\[
	h_p(w_1^p, w_2^p) = \frac{1}{\sqrt{\pi}} \int_{-\infty}^{\infty} e^{-t_{new}^2 } h_{1} h_{2} \text{d}t_{new}
\]
is then evaluated analytically at each $w_1^p$ and $w_2^p$ interpolation point. Finally,
another 2D FFT is performed to derive the coefficients of $h_p$. At the root level, the 
constant function (with no $t$-variables) gives the result of the integral.
In the implementation, as we use unified formulas for both the boundary nodes and interior nodes, 
the two functions $leftbdry$ and $rightbdry$ become unnecessary, see Appendix for details. 
Except for the detailed implementations in the functions $leafnode$, 
$root$, and $interiornode$, the recursively implemented Matlab algorithm for the upward pass is 
identical in structure as the presented Algorithm $1$ for the tri-diagonal case, we therefore 
skip the pseudo-code.

The algorithm complexity can be computed as follows. In both the upward pass and downward pass,
constant numbers of operations and storage are required for each tree node, the overall algorithm
complexity and memory requirement are therefore both asymptotically optimal $O(N)$ for the
$N$-dimensional integration problem.

\subsection{Preliminary Numerical Results}
\label{sec:num2}
We present some preliminary results to demonstrate the accuracy and efficiency of the numerical algorithm
for the exponential case. The $N$ $z$-location points are randomly chosen in $[0,1]$ and sorted. 
A uniform tree is then generated by recursively subdividing the $z$-locations and corresponding 
$x$-variables, and the same settings of $\ba$ and $\bb$ are used as in the tridiagonal case.
We first study the accuracy of the algorithm. For $N=4$, we compute a reference solution
$\phi= 9.63128791560604001$ using Mathematica, with an estimated error $5.99\text{e-}8$. 
For $N=8$, Mathematica returns the result $\phi=1.16750673314578e+02 $ with an estimated
error $0.064$. For $N>8$, direct computation using Mathematica becomes impossible.
In Table \ref{tab:expAccuracy}, we show the Matlab results for different dimensions $N$
when $2M$ Fourier series terms are used in the approximation. The error tolerance for the
NUFFT solver is set to $1\text{e-}12$. For all cases, our results converge when $M$ increases. 
For both $N=4$ and $N=8$, our converged results match those from Mathematica within the 
estimated error from Mathematica.

In the current implementation, as the exponential case involves operations on a 3-variables 
function $h(w_1^p, w_2^p, t_{new})$ for each child when forming the parent's 
Fourier series expansion, while both the storage and operations for the tridiagonal case can be 
compressed so one only works on 2-variable functions (variables $\{t_l, t_m\}$ for child $1$ 
and $\{t_m, t_r\}$ for child $2$), the exponential solver therefore requires 
more operations and memory than the tridiagonal case. We tested our code on a 
desktop with 16GB memory and Intel Xeon CPU E3-1225 v6 @3.30GHz. 
For $N=4$ and $M=512$, More than $8G$ memory is already required, hence no result 
for $M=512$ is reported. 
\begin{table}[htbp]
\centering
\begin{tabular}{|c|c|c|c|}
\hline
$N$ & 4 & 8 & 16  \\
\hline
$M$=16  & 9.646301617204299 & 118.8260790816760  & 21594.43676761628   \\
$M$=32  & 9.631244805483258 & 116.7475848966488  & 17592.18271523017   \\
$M$=64  & 9.631287915305332 & 116.7505122381643  & 17591.75082916860   \\
$M$=128 & 9.631287915311097 & 116.7505122544810  & 17591.75095515863  \\
$M$=256 & 9.631287915311061 & 116.7505122544801  & 17591.75095515877  \\
\hline
\hline
$N$ & 32 & 64 & 128  \\
\hline
$M$=16  & 1131582930.741270 & 4.332761307147880e+18  & 7.074841023044070e+37   \\
$M$=32  & 550963842.9679267 & 1.046292247268069e+18  & 9.380354831605098e+36   \\
$M$=64  & 540456718.9698794 & 8.163524406713720e+17  & 3.432262767034514e+36   \\
$M$=128 & 540456737.4129881 & 8.163182314210313e+17  & 3.394537652388589e+36   \\
$M$=256 & 540456737.4129064 & 8.163182314206217e+17  & 3.394537652164628e+36   \\
\hline
\end{tabular}
\caption{Computed $\phi$ values for different dimensions and numbers of Fourier terms, exponential case.
\label{tab:expAccuracy}}
\centering
\end{table}

\vspace{0.1in}
{\noindent \bf Remark:} We explain the large errors when $M=16$ (and $M=32$) for large $N$ 
values. When the dimension of the problem increases, its condition number also increases
exponentially.
For each leaf node, if we assume the numerical solution has a relative error $\epsilon$ 
in each leaf node function $h_{leaf}$, in the worst case, the relative error for the 
$N$ dimensional integral can be approximated by $(1+\epsilon)^N-1$ as 
the $N$ leaf node functions will be ``multiplied" together in the upward pass to get
the final integral value. Clearly, the condition number of the analytical problem grows 
exponentially as $N$ increases. In our current implementation, we set the error tolerance
of the NUFFT solver to $10^{-12}$ relative error. Therefore, a very rough estimate for
the error when $N=128$, assuming $M$ is large enough so the leaf node function 
$h_{leaf}$ is resolved to machine precision, is given by 
$(1+10^{-12})^{128} \approx 1 + 10^{-10}$, i.e., at most $10$ digits are correct if the 
worst case happens. Our numerical results show that for the same $N$ value, all the converged 
results match at least in the first $10$ significant digits in Table~\ref{tab:expAccuracy}.

We demonstrate the efficiency of our algorithm by presenting the Matlab simulation time for different
$M$ and $N$ values, and the unit for the CPU time is in seconds. The current Matlab code has not been 
fully vectorized or parallelized, and significant performance improvement in the prefactor of the 
$O(N)$ algorithm is expected from a future optimized code. However, the numerical results 
in Table~\ref{tab:expEfficiency} using our existing code sufficiently and clearly show the asymptotic 
algorithm complexity: the CPU time grows approximately linearly as the dimension $N$ increases, 
and it increases by a factor of approximately $8$ as $M$ doubles.
\begin{table}[htbp]
\centering
\begin{tabular}{|c|c|c|c|c|c|}
\hline
$N$ & 4 & 8 & 16 & 32 & 64 \\
\hline
$M=32$ CPU time & 1.03 & 2.96 & 6.67 & 14.1 & 29.2 \\ 
\hline
$M=64$ CPU time & 8.06 & 23.4 & 53.8 & 115 & 238 \\ 
\hline
$M=128$ CPU time & 65.2 & 191 & 445 & 949 & 1965 \\ 
\hline
\hline
$N$ & 128 & 256 & 512 & 1024 & 2048 \\
\hline
$M=32$ CPU time & 59.4 & 119 & 244 & 471 & 951 \\ 
\hline
$M=64$ CPU time & 491 & 988 & 1924 & 3889 & 7766 \\ 
\hline
$M=128$ CPU time & 4049 & 8029 & 16068 & 32465 & 65073 \\ 
\hline
\end{tabular}
\caption{CPU time (in seconds) for different $M$ and $N$ values, exponential case.
\label{tab:expEfficiency}}
\centering
\end{table}

\section{Generalizations and Limitations}
\label{sec:generalize}
In both the tridiagonal and exponential cases, we present the algorithms for the case $H(\bx)=constant$.
For a general $H$ with low-rank properties, i.e., 
$$ H(\bx) = \sum_{p=1}^P \prod_{k=1}^N u_{p, k}(x_k),$$
as $P$ is a small number, we can evaluate the expectation of each $p$ term $\prod_{k=1}^N u_{p, k}(x_k)$
and then add up the results. As the $x$-variables are already separated in the representation, 
the downward decoupling process can be performed the same as that in the tridiagonal or exponential 
case. At the finest level, the leaf node's function $h_{leaf}$ becomes 
$$h_{leaf}(\bt_{leaf}) = \int_{a_k}^{b_k} u_{p,k}(x_k) e^{-\frac12 \alpha_k x_k^2} 
  e^{i (\bt_{leaf}^T D_{leaf}) x_k} \text{d}x_k. $$ 
Note that analytical formula is in general not available for $h_{leaf}$, a numerical scheme has to be 
developed to compute the Fourier coefficients of $h_{leaf}$. This is clearly
numerically feasible as the integral is one-dimensional and $h_{leaf}$ is effectively a single 
variable function.

Next we consider more general $A$ matrices. We restrict our attention to the symmetric positive 
definite $\mathcal{H}$-matrices, and discuss the required low-rank and low dimensional properties in order 
for our method to become asymptotically optimal $O(N)$. 
A minimal requirement from the algorithm is that the off-diagonal matrices should be low rank.
Consider a parent's matrix $A$ with such low rank off-diagonals and the corresponding 
$x$-variables,
\begin{equation*}
A = 
\begin{bmatrix}
A_{l,1} & A_{l,2}=\bV \Lambda \bU^T \\
A_{l,3}=\bU \Lambda \bV^T & A_{l,4}
\end{bmatrix}, 
\bx = 
\begin{bmatrix}
\bx_{l,1}\\
\bx_{l,2}
\end{bmatrix}, 
\end{equation*}
where the first index $l$ is the current level of the block matrices and point sets, 
and we assume $\Lambda$ is low rank, $rank(\Lambda)=P$. 
Then we can rewrite the quadratic term in the exponent of the integrand as
\begin{equation*}
\begin{array}{rcl}
\bx^T A \bx &=& \bx_{l,1}^T A_{l,1} \bx_{l,1} + {\color{blue} \bx_{l,1}^T \bV \Lambda \bU^T \bx_{l,2} 
  + \bx_{l,2}^T \bU \Lambda \bV^T \bx_{l,1} } + \bx_{l,2}^T A_{l,4} \bx_{l,2} \\
                &=& {\color{red} (B \bU^T \bx_{l,2} + B^{-T}\bV^T \bx_{l,1})^T \Lambda 
                     (B\bU^T\bx_{l,2} + B^{-T}\bV^T\bx_{l,1}) }  + \\
                & & \bx_{l,1}^T A_{l,1} \bx_{l,1}  - 
                         {\color{blue} \bx_{l,2}^T\bU B^T \Lambda B\bU^T \bx_{l,2} } + \\
                & & \bx_{l,2}^T A_{l,4} \bx_{l,2} - 
                         {\color{blue} \bx_{l,1}^T\bV B^{-1} \Lambda B^{-T}\bV^T \bx_{l,1} } \\
                &=& {\color{green} \bx_{l,1}^T (A_{l,1} - \bV B^{-1} \Lambda B^{-T} \bV^T) \bx_{l,1}} + \\
                & & {\color{green} \bx_{l,2}^T (A_{l,4} - \bU B^T \Lambda B \bU^T) \bx_{l,2} } + \\
                & & {\color{red} (B\bU^T \bx_{l,2} + B^{-T} \bV^T \bx_{l,1})^T \Lambda
                      (B \bU^T \bx_{l,2} + B^{-T} \bV^T \bx_{l,1}) }, 
\end{array}
\end{equation*}
where the first two {\color{green} green} terms are the child problems to be processed recursively at finer levels
after we use a number $P$ of $t$-variables to decouple the $\bx_{l,1}$ and $\bx_{l,2}$ variables using
Eq.~(\ref{eq:Fourier}). Clearly, the number of effective variables cannot be smaller than $P$ in this case.
There are several difficulties in this divide-and-conquer strategy. First, the 
$P \times P$ constant matrix $B$ should be chosen so that the resulting children's 
matrices are also symmetric positive definite. As the choice of $B$ is not unique, 
its computation is currently done numerically using numerical linear algebra tools, and we
are still searching for additional conditions so that we can have uniqueness in $B$ and 
{\it better} numerical stabilities in the algorithm. Second, consider a covariance matrix
of a general data set, compared with the original off-diagonal matrix blocks in $A_{l,1}$ and 
$A_{l,4}$, the numerical rank of the off-diagonal blocks of the new child matrices 
$A_{l,1} - \bV B^{-1} \Lambda B^{-T} \bV^T$ and 
$A_{l,4} - \bU B^T \Lambda B\bU^T$, may increase. In the worst case, the
new rank can be as high as the old rank plus $P$. When this happens, the number of $t$-variables
required will increase rapidly when decoupling the finer level problems, and the number of effective
variables also increases dramatically. Fortunately, for many problems of interest today, the
singular vectors $\bU$ and $\bV$ also have special structures. For example, when the off-diagonal
covariance function can be well-approximated by a low
degree polynomial expansion using the separation of variables, then the singular vectors are just the discretized
versions of these polynomials, therefore the rank of all the old and new off-diagonal matrix blocks 
cannot be higher than the number of the polynomial basis functions, and the number of 
effective variables is also bounded by this number. 
In numerical linear algebra language, this means that all the left (or right) singular vectors of
the off-diagonal blocks belong to the same low-dimensional subspace, so that the  
singular vectors of the new child matrices $A_{l,1} - \bV B^{-1} \Lambda B^{-T} \bV^T$ and 
$A_{l,4} - \bU B^T \Lambda B \bU^T$ can be represented by the same set of
basis vectors in the subspace. For problems with this property, our algorithm can be generalized,
by numerically finding the relations between the effective variables in the downward
pass, and finding the parent's function coefficients using its children's in the upward pass.  
The numerical complexity of the resulting algorithm remains asymptotically optimal $O(N)$.

However, our algorithm also suffers from several severe limitations due to the lack of effective
tools for high dimensional problems. The main limitation is the large prefactor in the $O(N)$ 
complexity, as the prefactor grows exponentially when the rank of the off-diagonal blocks
(and hence the number of effective variables) increases. We presented the results when 
the number of effective variables are no more than $2$ in this paper. When this number increases to 
$5 \sim 20$, it may still be possible to introduce the sparse grid ideas 
\cite{barthelmann2000high,gerstner1998numerical,nobile2008sparse,shen2010efficient} when integrating 
the multi-variable $h$-functions. When this number is more than $20$, as far as we know, no current 
techniques can analytically handle problems of this size. Also, notice that for the current 
numerical implementation of the exponential case, fast algorithms such as FFT and NUFFT have to be 
introduced or the computation will become very expensive. However, as far as we know, existing 
NUFFT tools are only available in $1$, $2$, and $3$ dimensions. Finally, as the condition number
of the problem increases exponentially as $N$ increases, it is important to have very accurate
representations of the $h$-functions for the hierarchical tree nodes so reasonable accurate 
results are possible in higher dimensions. We are currently studying possible 
strategies to overcome these hurdles, by studying smaller matrix blocks so the rank can be 
lower, and more promisingly, by coupling the Monte Carlo approach with our divide-and-conquer 
strategy \cite{genton2018hierarchical}. Results along these directions will be reported in the future.

\section{Conclusions}
\label{sec:conclusion}
The main contribution of this paper is an asympotically optimal $O(N)$ algorithm for evaluating the 
expectation of a function $H(\bX)$  
$$\phi(\ba, \bb; A) = \int_{\ba}^{\bb} H(\bx) f(\bx | A) \text{d}\bx, $$ 
where $f(\bx | A)$ is 
the truncated multi-variate normal distribution with zero mean for the $N$-dimensional 
random vector $\bX$, when the off-diagonal blocks of $A$ are ``low-rank" with ``low-dimensional" 
features and $H(\bx)$ is ``low-rank". In the algorithm, a downward pass is performed to obtain 
the relations between the parent's and children's effective variables, followed by an upward 
pass to construct the $h$-functions for each node on the hierarchical tree structure. 
The function at the tree root returns the desired expectation. Numerical results are presented 
to demonstrate the accuracy and efficiency of the algorithm. The generalizations and limitations 
of the new algorithm are also discussed, with possible strategies so the algorithm can be applied 
to a wider class of problems.

\section*{Acknowledgement} 
J. Huang was supported by the NSF grant DMS1821093, and the work was finished while he was visiting professors at the 
King Abdullah University of Science and Technology, National Center for Theoretical Sciences (NCTS) in Taiwan, 
Mathematical Center for Interdisciplinary Research of Soochow University, and Institute for Mathematical Sciences 
of the National University of Singapore.

\section*{Appendix}
We first present the detailed formulas for the Green's function $G_p(z,y)$ of a parent node $p$ and the functions $G_1(z,y)$ 
and $G_2(z,y)$ of $p$'s left child $1$ and right child $2$. These functions are defined as
$$G_p(z,y) = \left\{ \begin{array}{l} 
	coef^p \cdot g_r^p(z) \cdot g_l^p(y), \quad y<z, \\
	coef^p \cdot g_l^p(z) \cdot g_r^p(y), \quad y>z, 
\end{array} \right.
$$
$$G_1(z,y) = \left\{ \begin{array}{l} 
	coef^1 \cdot g_r^1(z) \cdot g_l^1(y), \quad y<z, \\
	coef^1 \cdot g_l^1(z) \cdot g_r^1(y), \quad y>z, 
\end{array} \right.
$$
$$G_2(z,y) = \left\{ \begin{array}{l} 
	coef^2 \cdot g_r^2(z) \cdot g_2^p(y), \quad y<z, \\
	coef^2 \cdot g_l^2(z) \cdot g_2^p(y), \quad y>z. 
\end{array} \right.
$$
We assume parent's $z$-locations satisfy $z \in [a,b]$. We choose $\zeta=c$ to separate the parent's locations, and 
the child intervals are therefore $[a,c]$ and $[c,b]$, respectively.

\vspace{0.1in}
{\noindent \bf Case 1: $p$ is the root node ($a=0$, $b=1$):} The functions are 
\begin{displaymath}
	\begin{array}{lll}
		g_l^p(z)=e^{z-1}, & g_r^p(z)=e^{1-z}, & \quad coef^p=\frac12; \\ 
		g_l^1(z)=e^{z-c}, & g_r^1(z)=\sinh(c-z), & \quad coef^1=1; \\ 
		g_l^2(z)=\sinh(z-c), & g_r^2(z)=e^{c-z}, & \quad coef^2=1; 
	\end{array}
\end{displaymath}
{\noindent \bf Case 2: $p$ is a left boundary node($a=0$):} The functions are 
\begin{displaymath}
	\begin{array}{lll}
		g_l^p(z)=e^{z-b}, & g_r^p(z)=\sinh(b-z), & \quad coef^p=1; \\ 
		g_l^1(z)=e^{z-c}, & g_r^1(z)=\sinh(c-z), & \quad coef^1=1; \\ 
		g_l^2(z)=\sinh(z-c), & g_r^2(z)=\sinh(b-z), & \quad coef^2=\frac{2 e^{b+c}}{e^{2 b}-e^{2 c}}; 
	\end{array}
\end{displaymath}
{\noindent \bf Case 3: $p$ is a right boundary node($b=1$):} The functions are 
\begin{displaymath}
	\begin{array}{lll}
		g_l^p(z)=\sinh(z-a), & g_r^p(z)=e^{a-z}, & \quad coef^p=1; \\ 
		g_l^1(z)=\sinh(z-a), & g_r^1(z)=\sinh(c-z), & \quad coef^1=\frac{2 e^{a+c}}{e^{2 c}-e^{2 a}}; \\ 
		g_l^2(z)=\sinh(z-c), & g_r^2(z)=e^{c-z}, & \quad coef^2=1;
	\end{array}
\end{displaymath}
{\noindent \bf Case 4: $p$ is an interior node:} The functions are 
\begin{displaymath}
	\begin{array}{lll}
		g_l^p(z)=\sinh(z-a), & g_r^p(z)=\sinh(b-z), & \quad coef^p=\frac{2 e^{a+b}}{e^{2 b}-e^{2 a}}; \\ 
		g_l^1(z)=\sinh(z-a), & g_r^1(z)=\sinh(c-z), & \quad coef^1=\frac{2 e^{a+c}}{e^{2 c}-e^{2 a}}; \\ 
		g_l^2(z)=\sinh(z-c), & g_r^2(z)=\sinh(b-z), & \quad coef^2=\frac{2 e^{b+c}}{e^{2 b}-e^{2 c}}; 
	\end{array}
\end{displaymath}

Next, we present the relations of the parent $p$'s two $w$-variables $w_1^p$ and $w_2^p$ with the left child  
$1$'s two $w$-variables \{$w_1^1$, $w_2^1$\} and right child $2$'s two $w$-variables \{$w_1^2$, $w_2^2$\}. We use
$t_{new}$ to represent the new $t$-variable introduced to divide the parent problem into two sub-problems of child $1$
and child $2$. We use a unified set of basis functions for each node on the hierarchical tree structure. 
For the parent node, the basis functions are $\{ \Phi_1^p=\cosh(z-c), \Phi_2^p=\frac{\sinh(z-c)}{b-a} \}$. The basis
functions for the left and right children are $\{ \Phi_1^1 = \cosh(z-p), \Phi_2^1=\frac{\sinh(z-p)}{c-a} \}$ and 
$\{ \Phi_1^2=\cosh(z-q), \Phi_2^2=\frac{\sinh(z-q)}{b-c} \}$, respectively, where $p$ and $q$ are 
either the interface $\zeta$ points when further subdividing the two child problems, or the mid-point of the child 
intervals when they become leaf nodes. 

\vspace{0.1in}
{\noindent \bf Case 1: $p$ is the root node ($a=0$, $b=1$):} Parent has no effective $w$-variables. 
\begin{displaymath}
	\begin{array}{l}
		w_1^1=\frac{t_{new} e^{p-c}}{\sqrt{2}}, \quad w_2^1=-\frac{t_{new} (a-c) e^{p-c}}{\sqrt{2}}; \\
		w_1^2=\frac{t_{new} e^{c-q}}{\sqrt{2}}, \quad w_2^2=-\frac{t_{new} (b-c) e^{c-q}}{\sqrt{2}}.
	\end{array}
\end{displaymath}
{\noindent \bf Case 2: $p$ is a left boundary node($a=0$):}
\begin{displaymath}
	\begin{array}{l}
		w_1^1=\frac{w_2^p \sinh (c-p)}{a-b}+\frac{e^p t_{new} \sqrt{e^{-2 c}-e^{-2 b}}}{\sqrt{2}}+w_1^p \cosh (c-p), \\
	w_2^1=(a-c) \left(\frac{w_2^p \cosh (c-p)}{a-b}+
		w_1^p \sinh (c-p)\right)-\frac{e^p t_{new} (a-c) \sqrt{e^{-2 c}-e^{-2 b}}}{\sqrt{2}}; \\
		w_1^2= \frac{w_2^p \sinh (c-q)}{a-b}+t_{new} \sqrt{\coth (b-c)-1} \sinh (b-q)+w_1^p \cosh (c-q),  \\
		w_2^2= \frac{(b-c) (w_1^p (b-a) \sinh (c-q)-w_2^p \cosh (c-q))}{a-b}+t_{new} (c-b) \sqrt{\coth (b-c)-1} \cosh (b-q).
	\end{array}
\end{displaymath}
{\noindent \bf Case 3: $p$ is a right boundary node($b=1$):} 
\begin{displaymath}
	\begin{array}{l}
		w_1^1= \frac{w_2^p \sinh (c-p)}{a-b}+t_{new} \sqrt{-\coth (a-c)-1} \sinh (p-a)+w_1^p \cosh (c-p), \\
		w_2^1=(a-c) \left(\frac{w_2^p \cosh (c-p)}{a-b}+w_1^p \sinh (c-p)\right)+t_{new} (c-a) \sqrt{-\coth (a-c)-1} \cosh (a-p); \\
		w_1^2=\frac{w_2^p \sinh (c-q)}{a-b}+\frac{e^{-q} t_{new} \sqrt{e^{2 c}-e^{2 a}}}{\sqrt{2}}+w_1^p \cosh (c-q),  \\
		w_2^2=\frac{e^{-q} t_{new} \sqrt{e^{2 c}-e^{2 a}} (c-b)}{\sqrt{2}}+\frac{(b-c) (w_1^p (b-a) \sinh (c-q)-w_2^p \cosh (c-q))}{a-b}.
	\end{array}
\end{displaymath}
{\noindent \bf Case 4: $p$ is an interior node:} 
\begin{displaymath}
	\begin{array}{l}
		w_1^1=t_{new} \sinh (p-a) \sqrt{\text{csch}(a-b) \text{csch}(a-c) \sinh (b-c)}+\frac{w_2^p \sinh (c-p)}{a-b}+w_1^p \cosh (c-p) , \\
		w_2^1= t_{new} (c-a) \cosh (a-p) \sqrt{\text{csch}(a-b) \text{csch}(a-c) \sinh (b-c)}+(a-c) \left(\frac{w_2^p \cosh (c-p)}{a-b}+w_1^p \sinh (c-p)\right); \\
		w_1^2= t_{new} \sinh (b-q) \sqrt{\text{csch}(a-b) \sinh (a-c) \text{csch}(b-c)}+\frac{w_2^p \sinh (c-q)}{a-b}+w_1^p \cosh (c-q), \\
		w_2^2=t_{new} (c-b) \cosh (b-q) \sqrt{\text{csch}(a-b) \sinh (a-c) \text{csch}(b-c)}+\frac{(b-c) (w_1^p (b-a) \sinh (c-q)-w_2^p \cosh (c-q))}{a-b}.
	\end{array}
\end{displaymath}
Mathematica files for computing these formulas are available.

\bibliographystyle{plain}
\bibliography{huangbib}

\end{document}